\documentclass[12pt]{article}
\usepackage{amsmath,amssymb}

\newtheorem{Assumption}{Assumption}[part]

\def\cal#1{\mathcal{#1}}

\def \F{I\!\!F}

\def \R{I\!\!R}

\def\Fc{{\cal F}}

\def\Dzw1#1{\frac{\partial^2 #1}{\partial z \partial w_1}}

\def\Dzb1#1{\frac{\partial^2 #1}{\partial z \partial b_1}}

\newcommand{\dproof}{\noindent {Proof.} \quad}
\newcommand{\fproof}{\hfill $\square$ \bigskip}

\newtheorem{definition}{Definition}[section]

\newtheorem{theorem}[definition]{Theorem}

\newtheorem{lemma}[definition]{Lemma}

\def\R{{\bf R}}

\def\1B{\text{1\!\!I}}

\def\R{{\bf R}}

\def\1B{\text{1\!\!I}}

\begin{document}

\title{Stochastic control of  mean-field SPDEs with jumps}
\author{ Roxana Dumitrescu\thanks{ Department of Mathematics, King's College London, United Kingdom, email: {\tt roxana.dumitrescu@kcl.ac.uk}} \and Bernt $\O$ksendal \thanks{ Department of Mathematics, University of Oslo, P.O. Box 1053 Blindern, N-0316 Oslo, Norway, email: {\tt oksendal@math.uio.no}} 
%\thanks{CEREMADE,
%Universit\'e Paris 9 Dauphine, CREST  and  INRIA Paris-Rocquencourt, email: {\tt roxana@ceremade.dauphine.fr}. The research leading to these results has received funding from the R\'egion Ile-de-France. }
\and
Agn\`es Sulem
\thanks{INRIA Paris,   Equipe-projet MathRisk, 3 rue Simone Iff, CS 42112, 75589 Paris Cedex 12, France, email: {\tt agnes.sulem@inria.fr}}}

\date{\today}

\maketitle

\begin{abstract}

We study the problem of optimal control for mean-field stochastic partial differential equations (stochastic evolution equations) driven by a Brownian motion and an independent Poisson random measure, in the case of \textit{partial information} control. One important novelty of our problem is represented by the introduction of  \textit{general mean-field} operators, acting on both the controlled state process and the control process.   We first formulate  a sufficient and  a necessary maximum principle for  this type of control. We then prove existence and uniqueness of the solution of such general forward and backward mean-field stochastic partial differential equations. We finally apply our results to find the explicit optimal control for an optimal harvesting problem.
\end{abstract}


\section{Introduction}\label{sec1}

%In this paper we consider an optimal control problem for a stochastic process $Y(t,x)=Y^{u}(t,x)$ which follows the following dynamics:
%
%\begin{align}
%dY(t,x)&=&AY(t,x)+f(t,x,Y(t,x), \textbf{F(}Y(t,x)\textbf{)}, u(t,x), \textbf{G(}u(t,x)\textbf{)})dt\nonumber \\&+& \sigma\left(t,x,Y(t,x),\textbf{F(}Y(t,x)\textbf{)},u(t,x),\textbf{G(}u(t,x)\textbf{)}\right)dW_t\nonumber \\&+&\int_\textbf{E}\theta(t,x,Y(t,x),\textbf{F(}Y(t,x)\textbf{)},u(t,x),\textbf{G(}u(t,x)\textbf{)},e)\Tilde{N}(dt,de),
%\end{align}
%
%where $\textbf{F}, \textbf{G}$ are operators acting on the random variables belonging to $\textbf{L}_2(\textbf{P})$. One important example is represented by the expectation operator, which is classicaly used in the mean-field problems.
%
%The boundary conditions are
%$$Y(0,x)=\xi(x), x \in D $$
%and
%$$Y(t,x)=\eta(t,x);\,\,\, (t,x) \in [0,T] \times \partial D.$$
%
%
%
%

\subsection{A motivating example}
As a motivation for the problem studied in this paper, we consider the following \textbf{optimal harvesting} problem: Suppose we model the density $Y(t,x)$ of a fish population in a lake $D$ at time $t$ and at point $x \in D$ by an equation of the form:
\begin{align}\label{eqintro}
dY(t,x)=&\mathbf{E}[Y(t,x)]b(t,x)dt+\dfrac{1}{2} \sum_{i=1}^{d}\dfrac{\partial^2}{\partial^2x_{i}}Y(t,x)dt+Y(t,x)\sigma(t,x)dW_t 
\nonumber \\&
+Y(t,x)\int_{\mathbf{R}^*} \theta(t,x,e) \Tilde{N}(dt,de).\nonumber \\
Y(0,x)=&y_0(x), \,\, x \in D,
\end{align}
where $D$ is a bounded domain in $\mathbf{R}^d$ and $y_0(x), b(t,x), \sigma(t,x), \theta(t,x,e)$ are  given bounded deterministic functions. Here $W_t$ is a Brownian motion and $\Tilde{N}(dt,de)=N(dt,de)-\nu(de)dt$ is an independent compensated Poisson random measure, respectively, on a filtered probability space $(\Omega, \mathcal{F}, \mathbb{F}=\{\mathcal{F}_t\}, P)$.

We may heuristically regard \eqref{eqintro} as a limit as $n \rightarrow \infty$ of a large population interacting system of the form
\begin{align}\label{eqintro1}
dy^{j,n}(t,x)&=\left[\dfrac{1}{n}\sum_{l=1}^n y^{l,n}(t,x)\right]b(t,x)dt+\dfrac{1}{2}\sum_{i=1}^{d}\dfrac{\partial^2}{\partial^2x_{i}}y^{j,n}(t,x)dt+y^{j,n}(t,x)\sigma(t,x)dW_t \nonumber \\&+y^{j,n}(t,x)\int_{\mathbf{R}^*} \theta(t,x,e) \Tilde{N}(dt,de),\,\,j=1,2,...,n \nonumber \\
y^{j,n}(t,x)(0,x)&=y_0(x),
\end{align}
where we have divided the whole lake into a grid of size $n$ and $y^{j,n}(t,x)$ represents the density in box $j$ of the grid.
Now suppose we introduce a harvesting-rate process $u(t,x)$. The density of the corresponding population $Y(t,x)=Y^u(t,x)$ is thus modeled by a controlled mean-field stochastic partial differential equation with jumps of the form:
\begin{align}
dY(t,x)=\mathbf{E}[Y(t,x)]b(t,x)dt+\dfrac{1}{2}\sum_{i=1}^d\dfrac{\partial^2}{\partial^2x_{i}}Y(t,x)dt+Y(t,x)\sigma(t,x)dW_t \nonumber \\
+Y(t,x)\int_{\mathbf{R}^*} \theta(t,x,e) \Tilde{N}(dt,de)-Y(t,x) u(t,x)dt.
\end{align}
%%where $u(t,x)$ represents the harvesting-rate process.

The \textit{performance functional} is assumed to be of the form
\begin{equation}\label{eqintro3}
J(u)=\mathbf{E} \left[\int_{0}^{T}\int_{D} \log(Y(t,x)u(t,x))dxdt+ \int_D \alpha(x)Y(T,x)dx\right].
\end{equation}
This may be regarded as the expected total logarithmic utility of the harvest up to time $T$ plus the value of the remaining population at time $T$.

The problem is thus to find $u^*$ such that
\begin{equation}\label{eqqintro}
J(u^*)=\sup_{u \in \mathcal{A}} J(u),
\end{equation}
where $\mathcal{A}$ represents the set of \textit{admissible} controls.
This process $u^*(t,x)$ is called an optimal harvesting rate.

This is an example of an optimal control problem of a mean-field stochastic reaction-diffusion equation. In the next sections, we will give necessary and sufficient conditions for optimality of a control in the case of \textit{partial information}, as well as results of existence and uniqueness for forward and backward mean-field stochastic partial differential equations with a \textit{general mean-field operator}. Finally, we apply our results in order to solve the optimal harvesting problem presented above.

To the best of our knowledge, the only paper that deals with optimal control of mean-field SPDEs is \cite{TM}.
Our paper extends \cite{TM} in four ways:
%\begin{itemize}
%\item[(i)] 
(i) we consider a more general mean-field operator;
%\item[(ii)] 
(ii) we introduce an additional general mean-field operator which acts  on the control process;
%\item[(iii)] 
(iii) we add jumps;
%\item[(iv)] 
(iv) we study the optimal control problem in the case of \textit{partial information}.
%\end{itemize}

The paper is organized as follows: in Section 2 we show the sufficient and necessary maximum principles and apply the results to the optimal harvesting example. In Section 3, we investigate the existence and the uniqueness of the solution of  mean-field SPDEs with jumps and general mean-field operator. In Section 4, we prove the existence and the uniqueness of mean-field backward SPDEs with jumps and general mean-field operator.

\section{Maximum principles for optimal control with partial information of general mean-field SPDEs with jumps}

\subsection{Framework and formulation of the optimal control problem}
  Let $(\Omega, \Fc_, \F = \{\Fc_t\}_{0 \leq t \leq T}, P)$ be a filtered probability space. 
Let  $W$ be   a  one-dimensional Brownian motion.  Let ${\bf E}:= \R^*$ and ${\cal B} ({\bf E})$ be 
its Borel filtration. Suppose that it is 
%$({\bf E}, {\cal B} ({\bf E}))$  be a measurable space 
equipped with a $\sigma$-finite positive measure $\nu$, satisfying $\int_{\textbf{E}} |e|^2 \nu(de)<\infty$ and
 let  $N(dt,de)$ be a independent Poisson random measure with compensator $\nu(de)dt$.
 We denote by $\tilde N(dt,de)$ its compensated process, defined as $\Tilde{N}(dt,de)=N(dt,de)-\nu(de)dt$. For simplicity, we consider $d=1.$

We introduce the following notation:
\begin{itemize}
\item $\textbf{L}^2(\textbf{P})$:= the set of random variables $X$ such that $\mathbf{E}[|X|^2]<\infty.$
%\item $\textbf{L}^2(\textbf{R})$:= the set of measurable functions $k:(D, \mathcal{B}%(D)) \rightarrow (\mathbf{R}, \mathcal{B}(\mathbf{R}))$ with $\int_{D} k^2(x)dx< \infty$.
\item $\textbf{L}^2(\textbf{R})$:= the set of measurable functions $k:(\textbf{R}, \mathcal{B}(\textbf{R})) \rightarrow (\mathbf{R}, \mathcal{B}(\mathbf{R}))$ with $\int_{\textbf{R}} k^2(x)dx< \infty$.
\item $\textbf{H}^2$:= the set of real-valued predictable processes $Z(t,x)$ with $\mathbf{E}[\int_{0}^{T} \int_D Z^2(t,x)dxdt]<\infty,$ where $D$ a bounded domain in $\textbf{R}$.

\item $\textbf{L}^2_{\nu}$:= the set of measurable functions $l:(\textbf{E}, \mathcal{B}(\textbf{E})) \rightarrow (\mathbf{R}, \mathcal{B}(\mathbf{R}))$ such that  $ \| l \|^2_{\textbf{L}^2_{\nu}}:=\int_{\textbf{E}}l^2(e) \nu(de)<\infty$. 
The set $\textbf{L}^2_{\nu}$ is a Hilbert space equipped with the scalar product $<l,l'>_{\nu}:= \int_{\textbf{E}}l(e)l'(e) \nu(de)$ for all $l,l' \in \textbf{L}^2_{\nu} \times \textbf{L}^2_{\nu}$.
\item $\textbf{H}^2_{\nu}$: the set of predictable real-valued processes $k(t,x,\cdot)$ with $\mathbf{E}[\int_0^T \int_D \|k(t,x,\cdot) \|_{\textbf{L}^2_{\nu}}]<\infty$.
\end{itemize}
Assume that we are given a subfiltration 
$$\mathcal{E}_t \subseteq \mathcal{F}_t; \,\,\, t \in [0,T], $$
representing the information available to the controller at time $t$. For example, we could have 
$$\mathcal{E}_t=\mathcal{F}_{(t-\delta)^+}\,\,\,\, ( \delta>0 \text{ constant}) $$
meaning that the controller gets a $\textit{ delayed}$ information flow compared to $\mathcal{F}_t$.\\
Consider a controlled mean-field stochastic partial differential equation $Y(t,x)=Y^{u}(t,x)$ at $(t,x)$  of the following form
\begin{align}\label{dyn1}
d&Y(t,x)=\big[LY(t,x)+b(t,x,Y(t,x),\textbf{F(}Y(t,x)\textbf{)}, u(t,x), \textbf{G(}u(t,x)\textbf{)})\big]dt \nonumber\\ &+\sigma(t,x,Y(t,x),\textbf{F(}Y(t,x)\textbf{)},u(t,x), \textbf{G(}u(t,x)\textbf{)})dW_t\nonumber \\ 
&+  \int_{\mathbf{E}} \theta(t,x,Y(t,x),\textbf{F(}Y(t,x)\textbf{)}, u(t,x),\textbf{G(}u(t,x)\textbf{)},e) \Tilde{N}(dt,de);\,\,\, (t,x) \in (0,T) \times D.
\end{align}
with boundary conditions
\begin{align}
Y(0,x)=\xi(x);\,\,\, x \in D \label{initial} \\
%\end{align}
%\begin{align}
Y(t,x)=\eta(t,x);\,\,\, (t,x) \in (0,T) \times \partial D. \label{boundary}
\end{align}
We interpret $Y$ as a weak (variational) solution to \eqref{dyn1}, in the sense that for $\phi \in C_0^{\infty}(D),$
\begin{align}
<Y_t, \phi>_{\textbf{L}^2(D)} = <y_0, \phi>_{\textbf{L}^2(D)}+\int_0^t<Y_s, L^*\phi>ds+\int_0^t<b(s,Y_s), \phi>_{\textbf{L}^2(D)}ds+\nonumber \\ \int_0^t<\sigma(s,Y_s), \phi>_{\textbf{L}^2(D)}dW_s +\int_0^t \int_{\mathbf{E}}<\theta(s,Y_s,e),\phi>_{\textbf{L}^2(D)} \Tilde{N}(ds,de),
\end{align}
where $L^*$ corresponds to the adjoint operator of $L$ and $<\cdot,\cdot>$ represents the duality product between $W^{1,2}(D)$ and $W^{1,2}(D)^*$, with $W^{1,2}(D)$ the Sobolev space of order 1. Existence and the uniqueness of the solution are proved in Section 4. 

Under this framework the It\^{o} formula can be applied to such SPDEs. See e.g. Pardoux \cite{P}, Pr\'evot  and Rockner \cite{PR}.
%Before giving the main result of existence and uniqueness, we recall the definition of a solution in this context.
%Here $y_0 \in K:=L^2(D).$ Let $V:=W_0^{1,2}$ be the Sobolev space of order 1 with zero boundary condition. Then $Y$ is understood as a weak (variational) solution, in the sense that $Y \in L^2([0,T];V)$  and for $\phi \in C^{\infty}_0(D),$

Here, $dY(t,x)=d_tY(t,x)$ is the differential with respect to $t$ and $L$ is a bounded linear differential operator acting on $x$. The process $u(t,x,\omega)$ is our \textit{control} process, taking values in a  set $\textbf{A} \subset \mathbf{R}$.  The functions $b:[0,T] \times \Omega \times D \times \mathbf{R} \times \textbf{A} \times \textbf{R} \mapsto \mathbf{R};\,\, (t, \omega, x, y, \bar{y},u, \bar{u}) \mapsto b(t, \omega, x, y, \bar{y},u, \bar{u})$, $\sigma:[0,T] \times \Omega \times D \times \mathbf{R} \times \textbf{A} \times \textbf{R} \mapsto \mathbf{R};\,\, (t, \omega, x, y, \bar{y},u, \bar{u}) \mapsto \sigma(t, \omega, x, y, \bar{y},u, \bar{u})$, $\theta:[0,T] \times \Omega \times D \times \mathbf{R} \times \textbf{A} \times \textbf{R} \times \textbf{E} \mapsto \mathbf{R};\,\, (t, \omega, x, y, \bar{y},u, \bar{u},e) \mapsto \theta(t, \omega, x, y, \bar{y},u, \bar{u},e)$ are given $\mathcal{F}_t$-predictable processes for each $x, y,\bar{y}, \bar{u} \in \mathbf{R}$, $u \in \textbf{A}$, $e \in \textbf{E}$. We assume that  $b, \sigma, \theta$ are $C_{b}^1$ and have linear growth with respect to $y, \bar{y}, u, \bar{u}$. We denote by $\mathcal{A}_{\mathcal{E}}$  a given family of \textit{admissible} controls, contained in the set of $\mathcal{E}_t$-predictable  stochastic processes $u(t,x) \in \textbf{A}$ satisfying $\mathbf{E}[\int_0^T\int_D u^2(t,x)dxdt]<\infty$ and such that  \eqref{dyn1}-\eqref{initial}-\eqref{boundary} has a unique c\`adl\`ag solution $Y(t,x)$.

 In the above equation, $\textbf{F}, \textbf{G}:\textbf{L}^2(\textbf{P}) \mapsto \mathbf{\mathbf{R}}$ are Fr\'echet differentiable operators. One important example is represented by the expectation operator $\mathbf{E}[\cdot]$.

Let $f:[0,T] \times D \times \mathbf{R}_2 \times \textbf{A} \times \mathbf{R} \mapsto \mathbf{R}$ and $g:D \times \mathbf{R}_2 \mapsto \mathbf{R}$ be a given profit rate function and bequest rate function, respectively. Moreover, we suppose that
\begin{align*}
\mathbf{E}\left[ \int_0^T \left( \int_D |f(t,x,Y(t,x),\textbf{F(}Y(t,x)\textbf{)},u(t,x), \textbf{G(}u(t,x)\textbf{)})|dx \right)dt  \right. \nonumber \\ \left. +\int_D |g(x,Y(T,x),\textbf{F(}Y(t,x)\textbf{)}|dx\right] < \infty,
\end{align*}
where $f(t,x,y,\bar{y}, u, \bar{u}),$ $g(x,y, \bar{y})$ are of class $C^1_b$ with respect to $(y,\bar{y},u, \bar{u})$ and continuous w.r.t to $t$. $\mathbf{E}$ denotes the expectation with respect to $P$.

For each $u \in \mathcal{A}_{\mathcal{E}}$, we define the \textit{performance functional} $J(u)$ by
\begin{align}
J(u)=\mathbf{E}\left[ \int_0^T \left( \int_D f(t,x, Y(t,x), \textbf{F(}Y(t,x)\textbf{)}, u(t,x),\textbf{G(}u(t,x)\textbf{)})dx \right)dt \right. \nonumber \\+\left. \int_D g(x,Y(T,x), \textbf{F(}Y(t,x)\textbf{)}|dx\right].
\end{align}
We aim  to maximize $J(u)$ over all $u \in \mathcal{A}_{\mathcal{E}}$ and our problem is the following:

 Find $u^* \in \mathcal{A}_{\mathcal{E}}$ such that
\begin{equation}\label{formprob}
\sup_{u \in \mathcal{A}} J(u)=J(u^*).
\end{equation}
Such a process $u^*$ is called an optimal control (if it exists), and the number
$
J=J(u^*)
$
is the \textit{value} of this problem.

\subsection{Sufficient maximum principle for partial information optimal control for mean-field SPDEs with jumps}

In this section, we prove necessary and sufficient maximum principles for optimal control with  $\textit{ partial information}$ in the case of a process described by a \textit{mean-field} stochastic \textit{partial} differential equation (in short MFSPDE) driven by a Brownian motion $(W)$ and a Poisson random measure $\Tilde{N}$.  The drift and the diffusion coefficients as well as the performance functional depend not only on the state and the control but also on the distribution of the state process, and also on the one of the control.
%Let $\mathcal{A}_{\mathcal{E}}$ denote a given family of controls, contained in the set of $\mathcal{E}_t$-predictable controls $u(\cdot)$ such that    has a variational solution. If $u \in \mathcal{A}_{\mathcal{E}},$ we call $u$ an admissible control.

%Our aim in this subsection is to prove a verification theorem for \textit{partial information} optimal control of a process described by a \textit{mean-field} stochastic \textit{partial} differential equation (in short MFSPDE) driven by a Brownian motion $(W)$ and a Poisson random measure $\Tilde{N}$. The verification theorem takes the form of  \textit{a sufficient stochastic maximum principle} (SMP).\\
Define the \textit{Hamiltonian} $H:[0,T] \times D  \times \mathbf{R}_2 \times \mathbf{A} \times \mathbf{R}_3 \times \textbf{L}^2_\nu \mapsto \mathbf{R}$ as follows:
\begin{align}
H(t,x,y, \bar{y},u,\bar{u},p,q,\gamma)=f(t,x,y, \bar{y},u, \bar{u})+b(t,x,y, \bar{y}, u, \bar{u})p+\sigma(t,x, y, \bar{y},u, \bar{u})q \nonumber \\+\int_\textbf{E} \theta(t,x,y, \bar{y},u, \bar{u},e)\gamma(e) \nu(e).
\end{align}

In our case, since the state process and the cost functional are of mean-field type, it turns out that the  adjoint equation will be a \textit{mean-field backward SPDE}, denoted  in the sequel MFBSPDEs.

We now introduce the \textit{adjoint} operator of the operator $L$, denoted by $L^*$, which satisfies
\begin{equation}\label{opadj}
(L^*\phi,\psi)=(\phi,L \psi), \,\,\, \text{ for all } \phi, \psi \in C_0^\infty(\mathbf{R}),
\end{equation}
where 
$$<\phi_1,\phi_2>_{\textbf{L}^2(\textbf{R})}:=(\phi_1, \phi_2)=\int_{\mathbf{R}} \phi_1(x) \phi_2(x)dx $$  is the inner product in $ \mathbf{L}^2(\mathbf{R}).$

For $u \in \mathcal{A}_{\mathcal{E}}$, we consider the following mean field backward stochastic partial differential equation (the \textit{adjoint} equation) in the three unknown processes $p(t,x) \in \mathbf{R}, q(t,x) \in \mathbf{R}, \gamma(t,x,\cdot) \in \textbf{L}^2_\nu;$ called the \textit{adjoint processes}:
\begin{align}
&dp(t,x)=- \left[L^*p(t,x)+\dfrac{\partial H}{\partial y}(t,x,Y(t,x),\textbf{F(}Y(t,x)\textbf{)},u(t,x),\textbf{G(}u(t,x)\textbf{)},p(t,x),q(t,x),\gamma(t,x,\cdot))\right]dt  \nonumber  \\
&-\mathbf{E}\left[\dfrac{\partial H}{\partial \bar{y}}(t,x,Y(t,x),\textbf{F(}Y(t,x)\textbf{)},u(t,x),\textbf{G(}u(t,x)\textbf{)},p(t,x),q(t,x),\gamma(t,x,\cdot))\right]\nabla\textbf{F(}Y(t,x)\textbf{)}dt \nonumber  \\
&+q(t,x)dW_t+\int_\mathbf{E}\gamma(t,x,e)\Tilde{N}(dt,de);\,\, (t,x) \in (0,T) \times D. 
\label{adjdyn} \\
&p(T,x)=\frac{\partial g}{\partial y}(x,Y(T,x), \textbf{F(}Y(T,x)\textbf{)})+\mathbf{E}\left[\frac{\partial g}{\partial \bar{y}}(x,Y(T,x), \textbf{F(}Y(T,x)\textbf{)}\right] \nabla\textbf{F(}Y(T,x)\textbf{)}; \,\, x \in D \label{adjterm} \\
&p(t,x)=0; \,\, (t,x) \in (0,T) \times \partial D. \label{adjbound}
\end{align}
Note that \eqref{adjdyn}  is equivalent to

\begin{align*}
dp(t,x)&=- \left[L^*p(t,x)+\dfrac{\partial f}{\partial y}(t,x,Y(t,x),\textbf{F(}Y(t,x)\textbf{)},u(t,x), \textbf{G(}u(t,x)\textbf{)})\right.  \nonumber\\ 
& \left. +\dfrac{\partial b}{\partial y}(t,x,Y(t,x),\textbf{F(}Y(t,x)\textbf{)},u(t,x), \textbf{G(}u(t,x)\textbf{)})p(t,x)\right]dt \nonumber \\
&-\left[\dfrac{\partial \sigma}{\partial y}(t,x,Y(t,x),\textbf{F(}Y(t,x)\textbf{)},u(t,x), \textbf{G(}u(t,x)\textbf{)})q(t,x)\right. \nonumber \\ 
& \left. +\int_{\mathbf{E}}\dfrac{\partial \theta}{\partial y}(t,x,Y(t,x),\textbf{F(}Y(t,x)\textbf{)},u(t,x), \textbf{G(}u(t,x)\textbf{)},e)\gamma(t,x,e)\nu(de)   \right]dt  \nonumber \\
&-\mathbf{E}\left[\dfrac{\partial f}{\partial \bar{y}}(t,x,Y(t,x),\textbf{F(}Y(t,x)\textbf{)},u(t,x), \textbf{G(}u(t,x)\textbf{)}) \right. \nonumber\\
\end{align*}
\begin{align} 
& \left. +\dfrac{\partial b}{\partial \bar{y}}(t,x,Y(t,x),\textbf{F(}Y(t,x)\textbf{)},u(t,x),\textbf{G(}u(t,x)\textbf{)})p(t,x)\right]\nabla\textbf{F(}Y(t,x)\textbf{)}dt  \nonumber \\
&-\mathbf{E}\left[\dfrac{\partial \sigma}{\partial \bar{y}}(t,x,Y(t,x),\textbf{F(}Y(t,x)\textbf{)},u(t,x), \textbf{G(}u(t,x)\textbf{)})q(t,x)\right]\nabla\textbf{F(}Y(t,x)\textbf{)}dt  \nonumber \\
&-\mathbf{E}\left[\int_{\mathbf{E}}\dfrac{\partial \theta}{\partial \bar{y}}(t,x,Y(t,x),\textbf{F(}Y(t,x)\textbf{)},u(t,x),\textbf{G(}u(t,x)\textbf{)},e)\gamma(t,x,e)\nu(de)\right]\nabla\textbf{F(}Y(t,x)\textbf{)}dt \nonumber \\
&+q(t,x)dW_t+\int_\mathbf{E}\gamma(t,x,e)\Tilde{N}(dt,de), \,\, x \in D. \label{adjoint}
\end{align}

%\vspace{5mm}

We now show the sufficient maximum principle.

\begin{theorem}[Sufficient Maximum Principle for mean-field SPDEs with jumps]\label{SuffMax}
Let $\hat{u} \in \mathcal{A}_{\mathcal{E}}$ with corresponding solution $\hat{Y}(t,x)$ and suppose that $\hat{p}(t,x), \hat{q}(t,x)$ and $\hat{\gamma}(t,x,\cdot)$ is a solution of the adjoint MFBSPDE \eqref{adjdyn}-\eqref{adjterm}-\eqref{adjbound}. Assume the following hold:

\begin{itemize}
\item[(i)] The maps $Y \mapsto g(x,Y, \textbf{F}(Y))$ and 
\begin{align}\label{concav}
(Y,u) \mapsto H(Y,u):=H(t,x,Y,\textbf{F}(Y), u, \textbf{G}(u),\hat{p}(t,x), \hat{q}(t,x), \hat{\gamma}(t,x,\cdot))
\end{align}
are concave functions with respect to $Y$ and $(Y,u)$, respectively, for all $(t,x) \in [0,T] \times \bar{D}$.
\item[(ii)](The maximum condition)
\begin{align}\label{maxcond}
\mathbf{E}\left[ H(t,x,\hat{Y}(t,x), \textbf{F}(\hat{Y}(t,x)), \hat{u}(t,x), \textbf{G}(\hat{u}(t,x)), \hat{p}(t,x),\hat{q}(t,x),\hat{\gamma}(t,x, \cdot)) | \mathcal{E}_t \right]= \nonumber \\
ess \sup_{v \in \mathcal{A}_{\mathcal{E}}} \mathbf{E}\left[H(t,x,\hat{Y}(t,x),\textbf{F}(\hat{Y}(t,x)),v(t,x),\textbf{G}(v(t,x)),\hat{p}(t,x),\hat{q}(t,x),\hat{\gamma}(t,x, \cdot))| \mathcal{E}_t\right] \text{   a.s}.
\end{align}
for all $t\in [0,T]$ and $x \in \bar{D}$.

%\item[(iii)] For all $u \in \mathcal{A}$ we have
%
%\begin{align}
%\mathbf{E} \left[\int_D \int_0^T (Y(x,t)-\hat{Y}(t,x))^2 \left[ \hat{q}(t,x)^2+\int_\mathbf{E}\hat{\gamma}(t,x,e)^2\nu(de)  \right]  \right]<\infty
%\end{align}
%%and

%\begin{equation}
%\mathbf{E}\left[\int_D\int_0^T \hat{p}(t,x)^2\sigma(t,x,Y(t,x),\textbf{F(}Y(t,x)\textbf{)},u(t,x))^2 dtdx\right]+
%\end{equation}
%
%$$\mathbf{E}\left[\int_D\int_0^T\int_\mathbf{E}\theta(t,x,Y(t,x),\textbf{F(}Y(t,x)\textbf{)},u(t,x),e)^2\nu(de)dtdx\right]<\infty.$$

\end{itemize}
Then $\hat{u}(t)$ is an optimal control for the random jump field control problem \eqref{formprob}.
\end{theorem}

\dproof
Define a sequence of stopping times $\tau_n; \,\, n=1,2,....$ as follows:
\begin{align*}
 \tau_n:= \inf \{t>0; \max \{  & \|\hat{p}(t)\|_{\textbf{L}^2(D)}, \|\hat{q}(t)\|_{\textbf{L}^2(D)}, \|\hat{\gamma}(t)\|_{\textbf{L}^2(D\times \textbf{E})}, \|\sigma(t)-\hat{\sigma}(t)\|_{\textbf{L}^2(D)}, \|\theta(t)-\hat{\theta}(t)\|_{\textbf{L}^2(D\times \textbf{E})},
  \nonumber \\ &\|Y(t)-\hat{Y}(t)\|_{\textbf{L}^2(D)} \geq n \} \wedge T.
\end{align*}
Then  $\tau_n \rightarrow T$ as $n \rightarrow \infty$ and 
\begin{align*}
&\textbf{E} \left[\int_0^{\tau_n} \left(\int_{D}\hat{p}(t,x)(\sigma(t,x)-\hat{\sigma}(t,x))dx \right)dW_t+\int_0^{\tau_n} \int_{\textbf{E}}\left(\int_{D}(\theta(t,x,e)-\hat{\theta}(t,x,e))dx \right) \Tilde{N}(dt,de) \right] \\
%\end{align}
%\begin{align*}
& = \textbf{E} \left[\int_0^{\tau_n} \left(\int_{D}(Y(t,x)-\hat{Y}(t,x))\hat{q}(t,x)dx\right)dW_t \right.  \\ 
 & \quad +\left. \int_0^{\tau_n}  \int_{\textbf{E}}\left(\int_{D}(Y(t,x)-\hat{Y}(t,x))\hat{\gamma}(t,x,e)dx\right)\Tilde{N}(dt,de) \right]=0 \,\,\, \text{ for all } n.
\end{align*}
%Define a sequence of stopping times $\tau_n;\,\, n=1,2,...,$ as follows
%
%\begin{align}
%\tau_n= \inf \{ t>0;\,\, \max\{ |Y(x,t)-\hat{Y}(t,x)|,\, |\hat{q}(t,x)|
%\end{align}
Let us fix $u \in \mathcal{A}_{\mathcal{E}}$ and let $Y(t,x)=Y^{u}(t,x)$ be the associated solution of \eqref{dyn1}. Define:
\begin{equation}
\begin{cases}
\hat{f}:=f(t,x,\hat{Y}(t,x), \textbf{F}(\hat{Y}(t,x)), \hat{u}(t,x), \textbf{G}(\hat{u}(t,x)));\,\,\, f:=f(t,x,Y(t,x), \textbf{F(}Y(t,x)\textbf{)}, u(t,x), \textbf{G}(u(t,x)));\\
\hat{g}:=g(x,\hat{Y}(T,x), \textbf{F}(\hat{Y}(T,x)));\,\,\,g:=g(x,Y(T,x), \textbf{F(}Y(T,x)\textbf{)});
\end{cases}
\end{equation}
and
$$
\begin{cases}
\hat{b}:=b(t,x,\hat{Y}(t,x), \textbf{F}(\hat{Y}(t,x)), \hat{u}(t,x),\textbf{G}(\hat{u}(t,x)));\,\,\, b:=b(t,x,Y(t,x), \textbf{F(}Y(t,x)\textbf{)}, u(t,x),\textbf{G(}u(t,x)\textbf{)});\\
\hat{\sigma}:=\sigma(t,x,\hat{Y}(t,x), \textbf{F}(\hat{Y}(t,x)), \hat{u}(t,x),\textbf{G}(\hat{u}(t,x)));\,\,\, \sigma:=\sigma(t,x,Y(t,x), \textbf{F(}Y(t,x)\textbf{)}, u(t,x),\textbf{G(}u(t,x)\textbf{)});\\
\hat{\theta}:=\theta(t,x,\hat{Y}(t,x), \textbf{F}(\hat{Y}(t,x)), \hat{u}(t,x),\textbf{G}(\hat{u}(t,x)),e);\,\,\, \nonumber \\ \theta:=\theta(t,x,Y(t,x), \textbf{F(}Y(t,x)\textbf{)}, u(t,x), \textbf{G(}u(t,x)\textbf{)}, e).
\end{cases}
$$
We also set
$$
\begin{cases}
\hat{H}:=H(t,x,\hat{Y}(t,x), \textbf{F}(\hat{Y}(t,x)),\hat{u}(t,x),\textbf{G}(\hat{u}(t,x)),\hat{p}(t,x),\hat{q}(t,x),\hat{\gamma}(t,x, \cdot));\\
H:=H(t,x,Y(t,x), \textbf{F(}Y(t,x)\textbf{)},u(t,x),\textbf{G(}u(t,x)\textbf{)},\hat{p}(t,x),\hat{q}(t,x),\hat{\gamma}(t,x, \cdot)).
\end{cases}
$$
Using the above definitions and the definition of the performance functional $J$, we get that:
\begin{equation}\label{first}
J(u)-J(\hat{u})=\mathcal{J}_1+\mathcal{J}_2,
\end{equation}
where $\mathcal{J}_1:=\mathbf{E}[\int_0^T\int_D(f-\hat{f})dxdt]$ and  $\mathcal{J}_2:=\mathbf{E}[\int_D(g-\hat{g})dx].$\\

Now, let us notice the following relations:
$$
\begin{cases}
\hat{f}=\hat{H}-\hat{b}\hat{p}(t,x)-\hat{\sigma}\hat{q}(t,x)-\int_\mathbf{E}\hat{\theta}\hat{\gamma}(t,x,e)\nu(de);\\
\hat{f}=H-b\hat{p}(t,x)-\sigma\hat{q}(t,x)-\int_\mathbf{E}\theta\hat{\gamma}(t,x,e)\nu(de),
\end{cases}
$$
which imply
\begin{equation}\label{second}
\mathcal{J}_1=\mathbf{E}\left[\int_0^T \int_D \left(H-\hat{H}-(b-\hat{b})\cdot \hat{p}-(\sigma-\hat{\sigma})\cdot \hat{q} -\int_\mathbf{E}(\theta-\hat{\theta})\cdot \hat{\gamma} \nu(de)\right)  \right].
\end{equation}
Fix $x \in D$. Since the map $Y \mapsto g(x,Y,F(Y))$ is concave for each $x \in \bar{D}$, we obtain:
$$
g-\hat{g} \leq \dfrac{\partial g}{\partial y}(x,\hat{Y}(T,x), \textbf{F}(\hat{Y}(T,x)))\Tilde{Y}(T,x)+\dfrac{\partial g}{\partial \bar{y}}(x,\hat{Y}(T,x), \textbf{F}(\hat{Y}(T,x)))<\nabla\textbf{F(}\hat{Y}(T,x)\textbf{)},\Tilde{Y}(T,x)>_{\textbf{L}^2(\textbf{P})},
$$
where
$$\Tilde{Y}(t,x)=Y(t,x)-\hat{Y}(t,x).$$
We thus obtain, by taking the expectation and applying the It\^o formula for jump-diffusion processes,
\begin{align}\label{third}
&\mathcal{J}_2 \leq \mathbf{E}\left[\int_D \left( \dfrac{\partial g}{\partial y}(x,\hat{Y}(T,x), \textbf{F}(\hat{Y}(T,x)))\Tilde{Y}(T,x) \right. \right. 
\nonumber \\
& \left. \left. \quad 
+\dfrac{\partial g}{\partial \bar{y}}(x,\hat{Y}(T,x), \textbf{F}(\hat{Y}(T,x)))<\nabla\textbf{F}(\hat{Y}(T,x)),\Tilde{Y}(T,x)>_{\textbf{L}^2(\textbf{P})}\right)dx\right] \nonumber \\
&= \mathbf{E}\left[\int_D <\hat{p}(T,x),\Tilde{Y}(T,x)> dx  \right]\nonumber \\
&  =\mathbf{E}\left[\int_D \left(\hat{p}(0,x)\cdot \Tilde{Y}(0,x)+\int_0^T \left( <\Tilde{Y}(t,x),d\hat{p}(t,x)>+ \hat{p}(t,x)d\Tilde{Y}(t,x)+(\sigma-\hat{\sigma})\hat{q}(t,x) \right)dt\right)dx\right] \nonumber \\
& \quad +\mathbf{E}\left[\int_D \left(\int_0^T\int_\mathbf{E}(\theta-\hat{\theta})\hat{\gamma}(t,x,e)N(dt,de) \right) dx\right] \nonumber \\
&=\mathbf{E} \left[ \int_D \int_0^T \hat{p}(t,x)\left(L\Tilde{Y}(t,x)+(b-\hat{b})\right) \right. \nonumber\\
& \quad +  \left. \Tilde{Y}(t,x)\left(-L^*\hat{p}(t,x)-\overline{\dfrac{\partial H}{\partial y}}-\mathbf{E}\left[\underline{\dfrac{\partial H}{\partial \bar{y}}}\right]<\nabla\textbf{F}(\hat{Y}(t,x)),\Tilde{Y}(t,x)>_{\textbf{L}^2(P)} \right) dt dx\right] \nonumber \\
& \quad +\mathbf{E}\left[\int_D  \int_0^T \left((\sigma-\hat{\sigma})\hat{q}(t,x)+\int_\mathbf{E} (\theta-\hat{\theta})\hat{\gamma}(t,x,e)\nu(de)\right)dtdx\right],
\end{align}
where
$$\overline{\dfrac{\partial H}{\partial y}}:= \dfrac{\partial H}{\partial y}(t,x,\hat{Y}(t,x),\textbf{F(}\hat{Y}(t,x)\textbf{)},\hat{u}(t,x),\textbf{G(}\hat{u}(t,x)\textbf{)},\hat{p}(t,x),\hat{q}(t,x),\hat{\gamma}(t,x,\cdot))$$
and
$$\underline{\dfrac{\partial H}{\partial \bar{y}}}:= \dfrac{\partial H}{\partial \bar{y}}(t,x,\hat{Y}(t,x),\textbf{F(}\hat{Y}(t,x)\textbf{)},\hat{u}(t,x),\textbf{G(}\hat{u}(t,x)\textbf{)},\hat{p}(t,x),\hat{q}(t,x),\hat{\gamma}(t,x,\cdot)).$$
From \eqref{first}, \eqref{second} and \eqref{third}, we derive
\begin{align*}
J(u)-J(\hat{u}) &\leq \mathbf{E}\left[\int_0^T \left(\int_D\hat{p}(t,x)L \widetilde{Y}(t,x)-\widetilde{Y}(t,x)L^*\hat{p}(t,x)dx\right)dt \right] \\
&+\mathbf{E}\left[\int_D \left(\int_0^T \left(H-\hat{H}-\overline{\dfrac{\partial H}{\partial y}}\cdot \Tilde{Y}(t,x)-\mathbf{E}\left[ \underline{\dfrac{\partial H}{\partial \bar{y}}} \right]<\nabla \textbf{F(}\hat{Y}(t,x)\textbf{)}, \Tilde{Y}(t,x)>_{\textbf{L}^2(\textbf{P})} \right)dt\right)dx\right]. 
\end{align*}
Since $\Tilde{Y}(t,x)=\hat{p}(t,x)=0$ for all $(t,x) \in [0,T] \times \partial D$, we obtain by an easy extension of  \eqref{opadj}  using Green's formula that
$$
\int_D \Tilde{Y}(t,x)L^*\hat{p}(t,x)dx=\int_D \hat{p}(t,x)L\Tilde{Y}(t,x),
$$ for all $t \in (0,T)$. 
We therefore get
\begin{align*}
J(u)-J(\hat{u}) \leq \mathbf{E}\left[\int_D \left(\int_0^T \left(H-\hat{H}-\overline{\dfrac{\partial H}{\partial y}}\cdot \Tilde{Y}(t,x)+\mathbf{E}\left[ \underline{\dfrac{\partial H}{\partial \bar{y}}} \right]<\nabla \textbf{F(}\hat{Y}(t,x)), \Tilde{Y}(t,x)>_{\textbf{L}^2(P)} \right)dt\right)dx\right].
\end{align*}
By the concavity assumption  \eqref{concav} we have
\begin{align*}
H-\hat{H} \leq \frac{\partial H}{ \partial y} (\hat{Y}, \textbf{F(}\hat{Y}\textbf{)},\hat{u}, \textbf{G(}\hat{u}\textbf{)})(Y-\hat{Y})+\frac{\partial H}{ \partial \bar{y}} (\hat{Y}, \textbf{F(}\hat{Y}\textbf{)},\hat{u}, \textbf{G(}\hat{u}\textbf{)})<\nabla \textbf{F(}\hat{Y}),(Y-\hat{Y})>_{\textbf{L}^2(\mathbf{P})} \nonumber \\+\frac{\partial H}{ \partial u} (\hat{Y}, \textbf{F(}\hat{Y}),\hat{u},\textbf{G(}\hat{u}))(u-\hat{u}) +\frac{\partial H}{ \partial \bar{u}} (\hat{Y}, \textbf{F(}\hat{Y}),\hat{u},\textbf{G(}\hat{u}))<\nabla \textbf{G(}\hat{u}),(u-\hat{u})>_{\textbf{L}^2(\mathbf{P})}.
\end{align*}
Combining the two above relations we get:
\begin{align}
J(u)-J(\hat{u}) 
 \leq &\mathbf{E}\left[\int_D \int_0^T \left( \frac{\partial H}{ \partial u} (\hat{Y}, \textbf{F(}\hat{Y}\textbf{)},\hat{u},\textbf{G(}\hat{u}))(u-\hat{u}) \right. \right.\nonumber \\
&\left.  \left. + \frac{\partial H}{ \partial \bar{u}} (\hat{Y}, \textbf{F(}\hat{Y}\textbf{)},\hat{u},\textbf{G(}\hat{u}))<\nabla \textbf{G(}\hat{u}),(u-\hat{u})>_{\textbf{L}^2(\textbf{P})}\right) dt dx \right]. \label{ineq}
\end{align}
By the maximum condition \eqref{maxcond} , we obtain:
\begin{align}\label{ineq2}
\mathbf{E}\left[\dfrac{ \partial H}{\partial u} (\hat{Y}, \textbf{F(}\hat{Y}\textbf{)},\hat{u},\textbf{G(}\hat{u}))| \mathcal{E}_t   \right](u-\hat{u})+\mathbf{E}\left[\dfrac{\partial H}{\partial \bar{u}}(\hat{Y}, \textbf{F(}\hat{Y}\textbf{)},\hat{u},\textbf{G(}\hat{u}))| \mathcal{E}_t   \right]<\nabla \textbf{G(}{\hat{u}}\textbf{)},u-\hat{u}>_{\textbf{L}^2(\mathbf{P})} \leq 0 \text{  a.s.},
\end{align}
$\text{for all } (t,x) \in [0,T] \times D$.
From \eqref{ineq} and \eqref{ineq2} we conclude that
\begin{align*}
J(u) \leq J(\hat{u}).
\end{align*}
By arbitrariness of $u$, we conclude that $\hat{u}$ is optimal.
\fproof

\subsection{A necessary-type maximum principle for partial information control of mean-field SPDEs with jumps}

 As in many applications the concavity condition may not hold, we prove a version of the maximum principle which does not need this assumption. 
Instead, we assume the following:

$\textbf{(A1)}$ For all $s \in [0,T)$ and all bounded $\mathcal{E}_s$-measurable random variables $\theta(\omega,x)$ the control $\beta$ defined by
$$\beta_t(\omega,x)=\theta(\omega,x)\chi_{(s,T]}(t); \,\, t \in [0,T], \,\, x \in D$$ is in $\mathcal{A}_{\mathcal{E}}$.

$\textbf{(A2)}$ For all $u, \beta \in \mathcal{U}$ where $\beta$ is bounded there exists $\delta>0$ such that the control

$$u(t)+y \beta(t);\,\, t \in [0,T]$$

belongs to $\mathcal{A}_{\mathcal{E}}$ for all $y \in (-\delta, \delta)$.\\

Let us give an  auxiliary lemma.

\begin{lemma}
Let $u \in \mathcal{A}_{\mathcal{E}}$ and $v \in \mathcal{A}_{\mathcal{E}}$. The derivative process
\begin{equation}\label{Equa}
\mathcal{Y}(t,x):=\lim_{z \mapsto 0^+} \dfrac{Y^{u+z \beta}(t,x)-Y^{u}(t,x)}{z}
\end{equation}
exists and belongs to $\textbf{L}^2(dx \times dt \times dP)$.
We then have that $\mathcal{Y}$ satisfies the following mean-field SPDE:
\begin{align*}
&d\mathcal{Y}(t,x) =L \mathcal{Y}(t,x)+ \left(\dfrac{\partial b}{\partial y}(t,x,Y^{u}(t,x), \textbf{F(}Y^{u}(t,x)\textbf{)}, \textbf{G(}u(t,x)\textbf{)}) \mathcal{Y}(t,x) \right. \nonumber \\
&+ \dfrac{\partial b}{\partial \bar{y}}(t,x,Y^{u}(t,x), \textbf{F(}Y^{u}(t,x)\textbf{)},u(t,x), \textbf{G(}u(t,x)\textbf{)})<\nabla \textbf{F(}Y^{u}(t,x)\textbf{)}, \mathcal{Y}(t,x)>_{\mathbf{L}^2(\textbf{P})}  \\
 &+\dfrac{\partial b}{\partial u}(t,x,Y^{u}(t,x), \textbf{F(}Y^{u}(t,x)\textbf{)},u(t,x),\textbf{G(}u(t,x)\textbf{)}) \beta(t,x) \\ 
 & \left.+\dfrac{\partial b}{\partial \bar{u}}(t,x,Y^{u}(t,x), \textbf{F(}Y^{u}(t,x)\textbf{)},u(t,x),\textbf{G(}u(t,x)\textbf{)})<\nabla \textbf{G(}u(t,x)\textbf{)}, \beta(t,x)>_{\mathbf{L}^2(\textbf{P})}\right)dt  \\
&+ \left(\dfrac{\partial \sigma}{\partial y}(t,x,Y^{u}(t,x), \textbf{F(}Y^{u}(t,x)\textbf{)},u(t,x),\textbf{G(}u(t,x)\textbf{)}) \mathcal{Y}(t,x) \right. \\
&+ \dfrac{\partial \sigma}{\partial \bar{y}}(t,x,Y^{u}(t,x), \textbf{F(}Y^{u}(t,x)\textbf{)},u(t,x),\textbf{G(}u(t,x)\textbf{)})<\nabla \textbf{F(}Y^{u}(t,x)\textbf{)}, \mathcal{Y}(t,x)>_{\mathbf{L}^2(\textbf{P})}  \\
 &+\dfrac{\partial \sigma}{\partial u}(t,x,Y^{u}(t,x), \textbf{F(}Y^{u}(t,x)\textbf{)},u(t,x), \textbf{G(}Y^{u}(t,x)\textbf{)}) \beta(t,x) \\ 
 & \left. + \dfrac{\partial \sigma}{\partial \bar{u}}(t,x,Y^{u}(t,x), \textbf{F(}Y^{u}(t,x)\textbf{)},u(t,x),\textbf{G(}u(t,x)\textbf{)})<\nabla \textbf{G(}u(t,x)\textbf{)}, \beta(t,x)>_{\mathbf{L}^2(\textbf{P})}\right)dW_t \\
 &+\int_\mathbf{E} \left(\dfrac{\partial \theta}{\partial y}(t,x,Y^{u}(t,x), \textbf{F(}Y^{u}(t,x)\textbf{)},u(t,x), \textbf{G(}u(t,x)\textbf{)},e) \mathcal{Y}(t,x) \right. 
\end{align*}
\begin{align*}
&+ \dfrac{\partial \theta}{\partial \bar{y}}(t,x,Y^{u}(t,x), \textbf{F(}Y^{u}(t,x)\textbf{)},u(t,x), \textbf{G(}u(t,x)\textbf{)},e)<\nabla \textbf{F(}Y^{u}(t,x)\textbf{)}, \mathcal{Y}(t,x)>_{\mathbf{L}^2(\textbf{P})} \nonumber\\
&+\dfrac{\partial \theta}{\partial u}(t,x,Y^{u}(t,x), \textbf{F(}Y^{u}(t,x)\textbf{)},u(t,x), \textbf{G(}u(t,x)\textbf{)},e) \beta(t,x) \nonumber \\
& \left. +\dfrac{\partial \theta}{\partial \bar{u}}(t,x,Y^{u}(t,x), \textbf{F(}Y^{u}(t,x)\textbf{)},u(t,x), \textbf{G(}u(t,x)\textbf{)},e) \beta(t,x))<\nabla \textbf{G(}u(t,x)\textbf{)}, \beta(t,x)>_{\mathbf{L}^2(\textbf{P})}\right.)\Tilde{N}(dt,de),  \\ 
&\mathcal{Y}(t,x)=0, \,\,\, (t,x) \in (0,T) \times \partial D; \\
&\mathcal{Y}(0,x)=0, \,\, x \in D.
\end{align*}
\end{lemma}
\dproof The result follows by applying the mean theorem. We omit the details.
\fproof

We now provide the necessary-type maximum principle for our optimal control problem for mean-field SPDEs.

\begin{theorem}[Necessary-type maximum principle for mean-field SPDEs with jumps]\label{necprin}
Let $\hat{u} \in \mathcal{A}_{\mathcal{E}}$ with corresponding solutions \eqref{dyn1}-\eqref{initial}-\eqref{boundary} and \eqref{adjdyn}-\eqref{adjterm}-\eqref{adjbound}. Assume that Assumptions \textbf{(A1)-\textbf{(A2)}} hold. Then the following are equivalent:\\

{\rm (i)} $ \dfrac{d}{dy} J(\hat{u}+y\beta)|_{y=0}=0  \text{ for all bounded } \beta \in \mathcal{A}_{\mathcal{E}}.$\\

%Suppose $\hat{u} \in \mathcal{A}_{\mathcal{E}}$ is optimal, i.e. $\sup_{u \in \mathcal{A}_{\mathcal{E}}} J(u)=J(\hat{u}).$ Let $\hat{Y}$, $(\hat{p},\hat{q},\hat{\gamma})$ be the corresponding solution of  $\eqref{dyn}$, respectively $\eqref{adjoint}$.  Then we get:\\

{\rm (ii)} $  \mathbf{E} \left[\nabla \hat{H}(t,x) |\mathcal{E}_t \right]= 0, \,\,\, \text{ for all } (t,x) \in [0,T] \times D \text{ a.s. },$\\

where
\begin{align*}
\nabla \hat{H}(t,x):=\dfrac{\partial H}{\partial u}(t,x, \hat{u}(t,x), \hat{Y}(t,x), \textbf{F(}\hat{Y}(t,x)\textbf{)}, \textbf{G(}\hat{u}(t,x)\textbf{)}, \hat{p}(t,x),\hat{q}(t,x),\hat{\gamma}(t,x, \cdot))\\+\mathbf{E}\left[\dfrac{\partial H}{\partial \bar{u}}(t,x, \hat{u}(t,x), \hat{Y}(t,x), \textbf{F(}\hat{Y}(t,x)\textbf{)}, \textbf{G(}\hat{u}(t,x)\textbf{)}, \hat{p}(t,x),\hat{q}(t,x),\hat{\gamma}(t,x, \cdot))\right] \nabla \textbf{G(}\hat{u}(t,x)\textbf{)},
\end{align*}
for all $(t,x) \in [0,T] \times D$.
\end{theorem}

\dproof 
The assumptions on the coefficients together with the mean theorem and relation \eqref{Equa} yield to:
\begin{align}
&\lim_{y \rightarrow 0} \dfrac{1}{y}\left(J(\hat{u}+y\beta)-J(\hat{u})  \right) = \mathbf{E}\left[ \int_0^T \int_{D} (\dfrac{\partial f}{\partial y}(t,x,\hat{Y}(t,x), \textbf{F(}\hat{Y}(t,x)\textbf{)}, \hat{u}(t,x), \textbf{G(}\hat{u}(t,x)\textbf{)}) \mathcal{Y}(t,x) \right. \nonumber\\ 
& \quad + \dfrac{\partial f}{\partial \overline{y}}(t,x,\hat{Y}(t,x), \textbf{F(}\hat{Y}(t,x)\textbf{)},\hat{u}(t,x),\textbf{G(}\hat{u}(t,x)\textbf{)})<\nabla \textbf{F(}\hat{Y}(t,x)\textbf{)}, \mathcal{Y}(t,x)>_{\textbf{L}^2(\textbf{P})} \nonumber \\
&  \quad+\dfrac{\partial f}{\partial u}(t,x,\hat{Y}(t,x), \textbf{F(}\hat{Y}(t,x)\textbf{)},\hat{u}(t,x),\textbf{G(}\hat{u}(t,x)\textbf{)}) \beta(t,x))dxdt]\nonumber \\
&  \quad+ \left. \dfrac{\partial f}{\partial \bar{u}}(t,x,\hat{Y}(t,x), \textbf{F(}\hat{Y}(t,x)\textbf{)},\hat{u}(t,x), \textbf{G(}\hat{u}(t,x)\textbf{)}) dxdt]<\nabla \textbf{G(}\hat{u}(t,x)\textbf{)}, \beta(t,x)>_{\textbf{L}^2(\textbf{P})} \right] \nonumber \\
& \quad +\mathbf{E}\left[\int_{D} (\dfrac{\partial g}{\partial y}(T,x,\hat{Y}(T,x), \textbf{F(}\hat{Y}(T,x)\textbf{)}) \mathcal{Y}(T,x) \right.  \nonumber \\
&  \quad +\left. \dfrac{\partial g}{\partial \overline{y}}(T,x,\hat{Y}(T,x), \textbf{F(}\hat{Y}(T,x)\textbf{)},\hat{u}(T,x))<\nabla \textbf{F(}\hat{Y}(T,x)\textbf{)}, \mathcal{Y}(T,x)>_{\textbf{L}^2(\textbf{P})}dx)\right].\label{eqa}
\end{align}
The definition of the Hamiltonian $H$  implies:
\begin{align}\label{eqq1}
&\dfrac{\partial f}{\partial y}(t,x,\hat{Y}(t,x), \textbf{F(}\hat{Y}(t,x)\textbf{)}, \hat{u}(t,x), \textbf{G(}\hat{u}(t,x)\textbf{)})= \nonumber \\ & = \dfrac{\partial \hat{H}}{\partial y}(t,x)-\frac{\partial \hat{b}}{\partial y}(t,x)\hat{p}(t,x)-\frac{\partial \hat{\sigma}}{\partial y}(t,x)\hat{q}(t,x)-\int_{\mathbf{E}} \frac{\partial \hat{\theta}}{\partial y}(t,x,e)\hat{\gamma}(t,x,e)\nu(de)
\end{align}
\begin{align}\label{eqq2}
&\dfrac{\partial f}{\partial \overline{y}}(t,x,\hat{Y}(t,x), \textbf{F(}\hat{Y}(t,x)\textbf{)},\hat{u}(t,x),\textbf{G(}\hat{u}(t,x)\textbf{)})= \nonumber \\ &= \frac{\partial \hat{H}}{\partial \overline{y}}(t,x)-\frac{\partial \hat{b}}{\partial \overline{y}}(t,x)\hat{p}(t,x)-\frac{\partial \hat{\sigma}}{\partial \overline{y}}(t,x)\hat{q}(t,x)-\int_{\mathbf{E}} \frac{\partial \hat{\theta}}{\partial \overline{y}}(t,x,e)\hat{\gamma}(t,x,e)\nu(de)
\end{align}
\begin{align}\label{eqq3}
& \dfrac{\partial f}{\partial u}(t,x,\hat{Y}(t,x), \textbf{F(}\hat{Y}(t,x)\textbf{)},\hat{u}(t,x),\textbf{G(}\hat{u}(t,x)\textbf{)})  \nonumber \\
&= \frac{\partial \hat{H}}{\partial u}(t,x)-\frac{\partial \hat{b}}{\partial u}(t,x)\hat{p}(t,x)-\frac{\partial \hat{\sigma}}{\partial u}(t,x)\hat{q}(t,x)-\int_{\mathbf{E}} \frac{\partial \hat{\theta}}{\partial u}(t,x,e)\hat{\gamma}(t,x,e)\nu(de)
\end{align}
\begin{align}\label{eqq4}
&\dfrac{\partial f}{\partial \bar{u}}(t,x,\hat{Y}(t,x), \textbf{F(}\hat{Y}(t,x)\textbf{)},\hat{u}(t,x), \textbf{G(}\hat{u}(t,x)\textbf{)}) \nonumber \\
&= \frac{\partial \hat{H}}{\partial \overline{u}}(t,x)-\frac{\partial \hat{b}}{\partial \overline{u}}(t,x)\hat{p}(t,x)-\frac{\partial \hat{\sigma}}{\partial \overline{u}}(t,x)\hat{q}(t,x)-\int_{\mathbf{E}} \frac{\partial \hat{\theta}}{\partial \overline{u}}(t,x,e)\hat{\gamma}(t,x,e)\nu(de)
\end{align}
Using \eqref{eqa}, \eqref{eqq1}, \eqref{eqq2}, \eqref{eqq3}, \eqref{eqq4} we derive:
\begin{align*}
&\lim_{y \rightarrow 0} \dfrac{1}{y}\left(J(\hat{u}+y\beta)-J(\hat{u})  \right) \\
& = \mathbf{E} \left[ \int_{0}^T \int_{D}  \left( \frac{\partial \hat{H}}{\partial y}(t,x)-\frac{\partial \hat{b}}{\partial y}(t,x)\hat{p}(t,x)-\frac{\partial \hat{\sigma}}{\partial y}(t,x)\hat{q}(t,x)-\int_{\mathbf{E}} \frac{\partial \hat{\theta}}{\partial y}(t,x,e)\hat{\gamma}(t,x,e)\nu(de)\right) \mathcal{Y}(t,x)dxdt\right] \\
&+\mathbf{E} [ \int_{0}^T \int_{D}  ( \frac{\partial \hat{H}}{\partial \overline{y}}(t,x)-\frac{\partial \hat{b}}{\partial \overline{y}}(t,x)\hat{p}(t,x)-\frac{\partial \hat{\sigma}}{\partial \overline{y}}(t,x)\hat{q}(t,x)-\int_{\mathbf{E}} \frac{\partial \hat{\theta}}{\partial \overline{y}}(t,x,e)\hat{\gamma}(t,x,e)\nu(de))   \\ 
& \qquad + <\!\!\nabla \textbf{F(}\hat{Y}(t,x)\textbf{)}, \mathcal{Y}(t,x)\!\!>_{\textbf{L}^2(\textbf{P})}dxdt ]  \\
&+\mathbf{E} \left[ \int_{0}^T \int_{D}  \left( \frac{\partial \hat{H}}{\partial u}(t,x)-\frac{\partial \hat{b}}{\partial u}(t,x)\hat{p}(t,x)-\frac{\partial \hat{\sigma}}{\partial u}(t,x)\hat{q}(t,x)-\int_{\mathbf{E}} \frac{\partial \hat{\theta}}{\partial u}(t,x,e)\hat{\gamma}(t,x,e)\nu(de)\right)\beta(t,x)dxdt\right] \nonumber \\
&+ \mathbf{E} [ \int_{0}^T \int_{D}  ( \frac{\partial \hat{H}}{\partial \overline{u}}(t,x)-\frac{\partial \hat{b}}{\partial \overline{u}}(t,x)\hat{p}(t,x)-\frac{\partial \hat{\sigma}}{\partial \overline{u}}(t,x)\hat{q}(t,x)-\int_{\mathbf{E}} \frac{\partial \hat{\theta}}{\partial \overline{u}}(t,x,e)\hat{\gamma}(t,x,e)\nu(de)) \\
& \qquad +  <\!\!\nabla \textbf{G(}\hat{u}(t,x)\textbf{)}, \beta(t,x)\!\!>dxdt] 
+  \mathbf{E} [  \int_{D} <\!\!\hat{p}(T,x), \mathcal{Y}(T,x)\!\!>dx].
\end{align*}
%By the above computations we get:
%$$\mathbf{E} \left[\int_0^T \int_D\frac{\partial \hat{H}}{\partial u}(t,x)\beta(t,x) dxdt\right]=0; \,\, \beta \in \mathcal{U}  \text{ is bounded }.$$
Applying It\^o formula to $<\hat{p}(T,x), \mathcal{Y}(T,x)>$ and using the dynamics of the adjoint equations, we finally get
\begin{equation}
\lim_{y \rightarrow 0} \dfrac{1}{y}\left(J(\hat{u}+y\beta)-J(\hat{u})  \right)=\mathbf{E} \left[ \int_{0}^T \int_{D} \mathbf{E}\left[<\nabla_u\hat{H}(t,x), \beta(t,x)>|\mathcal{E}_t \right]dxdt \right],
\end{equation}
where 
$$<\nabla_u \hat{H}(t,x), \beta(t,x)>=\dfrac{\partial H}{\partial u}(t,x) \beta(t,x)+\dfrac{\partial H}{\partial \bar{u}}(t,x) <\nabla \textbf{G(}\hat{u}\textbf{)},\beta(t,x)>_{\textbf{L}^2(\mathbf{P})}.$$
We conclude that 
$$
\lim_{y \rightarrow 0} \dfrac{1}{y}\left(J(\hat{u}+y \beta)-J(\hat{u})\right)=0
$$
if and only if
$$
\mathbf{E}\left[\int_0^T \int_D \mathbf{E}[<\!\!\nabla_{u} \hat{H}(t,x), \beta(t,x) \!\!>|\mathcal{E}_t]dxdt\right]=0.
$$
In particular this holds for all $\beta \in \mathcal{A}_{\mathcal{E}}$ which takes the form
$$
\beta(t,x)= \theta(\omega,x) \chi_{[s,T]}(t);\,\, t \in [0,T],
$$
for a fixed $s \in [0,T),$ where $\theta(\omega,x)$ is a bounded  $\mathcal{E}_s$-measurable random variable. We thus get that this is again equivalent to
$$\mathbf{E}\left[\int_s^T \int_D \mathbf{E}\left[ <\nabla_u \hat{H}(t,x), \theta>| \mathcal{E}_t\right]dx dt \right] = 0.$$
We now differentiate with respect to $s$ and derive that
$$\mathbf{E} \left[ \int_D \mathbf{E} \left[<\!\!\nabla_u\hat{H}(s,x), \theta\!\!> |\mathcal{E}_s\right] dx \right] = 0.$$
Since this holds for all bounded $\mathcal{E}_s$-measurable random variable $\theta$, we can easily conclude that
$$
\lim_{y \rightarrow 0} \dfrac{1}{y}\left(J(\hat{u}+y \beta)-J(\hat{u})\right)=0
$$
is equivalent to 
$$
\mathbf{E} \left[\dfrac{\partial \hat{H}}{\partial u}(t,x) | \mathcal{E}_t \right]+\mathbf{E} \left[\dfrac{\partial \hat{H}}{\partial {\bar{u}}}(t,x)\right]\nabla \textbf{G(}\hat{u}(t,x)\textbf{)}  = 0\,\, a.s., \text{ for all } (t,x) \in [0,T] \times D.
$$
\fproof
%which is equivalent to 
%
%$$\mathbf{E} \left[\dfrac{\partial \hat{H}}{\partial u}(t,x)+\mathbf{E}\left[\dfrac{\partial \hat{H}}{\partial \hat{u}}(t,x)\right] \nabla G(\hat{u}) | \mathcal{E}_t \right] = 0\,\, a.s., \text{ for all } (t,x) \in [0,T] \times D.$$

\subsection{Application to the optimal harvesting example}
We now return to the problem of optimal harvesting from a population in a lake $D$ stated in the motivating example. Thus we suppose the density $Y(t,x)$ of the population at time $t \in [0,T]$ and at point $x \in D$ is given by the stochastic reaction-diffusion equation \eqref{eqintro}, and the performance criterion is assumed to be as in \eqref{eqintro3}. For simplicity, we choose $d=1$ and $\mathcal{E}_t=\mathcal{F}_t$. In this case the Hamiltonian  gets the following form
$$%\begin{align}
H(t,x,y, \bar{y},u, \bar{u},p,q,\gamma)=\log(yu)+[b(t,x)\bar{y}-yu]p+\sigma(t,x)yq+\int_{\mathbf{R}^*} \theta(t,x,e)y\gamma(e)\nu(de), 
$$%\end{align}
and the adjoint BSDE becomes
\begin{align*}
&dp(t,x)  =[-\dfrac{1}{2} \dfrac{\partial^2}{\partial^2 x}p(t,x)+\dfrac{1}{Y(t,x)}+\sigma(t,x)q(t,x)+\int_{\mathbf{R}^*}\theta(t,x,e)\gamma(t,x,e)\nu(de)\\ 
 & \quad -u(t,x)p(t,x)-\mathbf{E}[b(t,x)p(t,x)] ]dt +q(t,x)dW_t+\int_{\mathbf{R}^*} \gamma(t,x,e)\Tilde{N}(dt,de)\\
&p(T,x)=\alpha(x), \,\, x \in D, \\
&p(t,x)=0, \,\,\, (t,x) \in (0,T) \times \partial D.
\end{align*}
We now apply the necessary maximum principle which implies the fact that if $u$ is an optimal control then  it satisfies the first order condition 
$$
u(t,x)=\dfrac{1}{Y(t,x)p(t,x)}.
$$
We summarize our results as follows:
\begin{theorem}
Assume that the conditions of Theorem \ref{necprin} hold.  Suppose a harvesting rate process $u(t,x)$ is  optimal  for the optimization problem \eqref{eqqintro}. Then
\begin{align}
u(t,x)=\dfrac{1}{Y(t,x) p(t,x)},
\end{align}
where $p(t,x)$ solves the MFBSPDE
\begin{align*}
&dp(t,x) =[-\dfrac{1}{2} \dfrac{\partial^2}{\partial^2 x}p(t,x)+\dfrac{1}{Y(t,x)}+\sigma(t,x)q(t,x)+\int_{\mathbf{R}^*}\theta(t,x,e)\gamma(t,x,e)\nu(de)-\mathbf{E}[b(t,x)p(t,x)] \\
%\right. \nonumber \\ \left. 
& \qquad \qquad -u(t,x)p(t,x) ]dt 
+q(t,x)dW_t+\int_{\mathbf{R}^*} \gamma(t,x,e)\Tilde{N}(dt,de)\\
&p(T,x)=\alpha(x), \,\, x \in D. \nonumber \\
&p(t,x)=0, \,\,\, (t,x) \in (0,T) \times \partial D. \nonumber
\end{align*}
\end{theorem}

\section{Existence and uniqueness of general forward  mean-field SPDEs with L\'evy noise}

We address here the problem of existence and uniqueness of the forward  mean-field SPDE \eqref{dyn1} with general mean-field operator, introduced in Section 2.
In order to do this, we first introduce the general framework.
Let $V, H$ be two separable Hilbert spaces such that $V$ is continously, densely imbedded in $H$. Identifying $H$ with its dual we have
$$V \subset H \approxeq H^* \subset V^* ,$$
where we have denoted by $V^*$ the topological dual of $V$. Let $L$ be a bounded linear operator from $V$ to $V^*$ satisfying the following coercivity hypothesis: There exist constants $\chi>0$ and 
$\zeta \geq 0$ such that
\begin{align}\label{coerc}
2<-Lu,u>+\chi|u|^2_{H} \geq \zeta || u||^2_{V}\,\,\, \text{ for all }  u \in V,
\end{align}
where $<Lu,u>=Lu(u)$ denotes the action of $Lu \in V^*$ on $u \in V$ and $|\cdot|_H$ (resp. $||\cdot||_{V}$) the norm associated to the Hilbert space $H$ (resp. $V$).

Let us introduce the notation adopted in this section.
\begin{itemize}
\item $\textbf{L}_{\nu}^2(H)$ is  the set of measurable functions $k:(\textbf{E}, \mathcal{B}(\textbf{E})) \mapsto (H,\mathcal{B}(H))$ such that\\ $|\!| k|\!|_{\textbf{L}_\nu^2(H)}:=\left(\int_{\textbf{E}}|k(e)|^2_{H}\nu(de)\right)^{\frac{1}{2}}<\infty$;
\item $\textbf{L}^2(\Omega, H)$ is the set of measurable functions $k:(\Omega,\mathcal{F}) \mapsto (H, \mathcal{B}(H))$ such that  $\mathbf{E}[|k|^2_{H}]<\infty$;

\item $\textbf{L}^2(\Omega, \textbf{L}^2_{\nu}(H))$ is the set of measurable functions $k:(\Omega,\mathcal{F}) \mapsto (\textbf{L}_{\nu}^2(H), \mathcal{B}(\textbf{L}_{\nu}^2(H)))$ such that $\textbf{E}[|\!|k|\!|^2_{\textbf{L}^2_{\nu}(H)}]<\infty$

\item $\textbf{L}^2(\Omega \times [0,T],H)$ (resp. $\textbf{L}_2(\Omega \times [0,T],V)$) is  the set of $\mathcal{F}_t$-adapted $H$-valued (resp. $V$-valued) processes $\Phi:\Omega \times [0,T]  \mapsto H$ (resp.V) such that $ \|\Phi \|^2_{\textbf{L}^2(\Omega \times [0,T],H)}:=\mathbf{E}[\int_0^T |\Phi(t)|^2_{H}dt] < \infty$ (resp. $ \|\Phi \|^2_{\textbf{L}^2(\Omega \times [0,T],V)}:=\mathbf{E}[\int_0^T \|\Phi(t) \|^2_{V}dt] < \infty$). 

\item $\textbf{L}^2(\Omega \times [0,T] \times \mathbf{E},H)$ is the set of  all the $\mathcal{P} \times \mathcal{B}(\mathbf{E})$-measurable $H$-valued maps $\theta: \Omega \times [0,T] \times \textbf{E}\mapsto H$ satisfying $ \|\theta \|_{\textbf{L}^2(\Omega \times [0,T] \times \mathbf{E},H)}:=\mathbf{E}[\int_0^T \int_{\textbf{E}}|\Phi(t,e)|^{2}_H\nu(de)dt]<\infty.$

\item $\textbf{S}^2(\Omega \times [0,T],H)$ denotes the set of $\mathcal{F}_t$-adapted $H$-valued cadlag processes $\Phi:\Omega \times [0,T] \mapsto H$ such that $ \|\Phi \|^2_{\textbf{S}^2(\Omega \times [0,T],H)}:=\mathbf{E}\left[ \sup_{0 \leq t \leq T} |\Phi(t)|^2_{H} \right]< \infty.$
\end{itemize}

\vspace{2mm}
The mean field SPDE under study is:
\begin{equation}\label{dyn}
dY_t=[LY_t+b(t,Y_t,\textbf{F(}Y_t)\textbf{)}]dt+\sigma(t,Y_t,\textbf{F(}Y_t\textbf{)})dW_t+  \int_{\mathbf{E}} \theta(t,Y_t,\textbf{F(}Y_t\textbf{)},e) \Tilde{N}(dt,de); \nonumber \\ \,\,\, (t,x) \in (0,T) \times D.
\end{equation}
We recall that this equation should be understood in the weak sense.
%where $b \in \textbf{L}_2([0,T],H),$  $\sigma \in \textbf{L}_2([0,T],H)$ and $\theta \in \textbf{L}_2([0,T] \times \textbf{E}, H ).$\\

Before giving the main result of this section, we make the following assumption on the coefficients $b, \sigma, \theta$ and the operator $\textbf{F}$ which appear in the above mean-field SPDE.
\begin{Assumption}\label{As1}
The maps $b:\Omega \times [0,T] \times H \times H \mapsto H$, $\sigma :\Omega \times [0,T] \times H \times H \mapsto H$ are $\mathcal{P} \times \mathcal{B}(H) \times \mathcal{B}(H)/\mathcal{B}(H)$-measurable. The map $\theta:\Omega \times [0,T]  \times H \times H \times \mathbf{E} \mapsto H$ is $\mathcal{P} \times \mathcal{B}(\mathbf{E}) \times \mathcal{B}(H) \times \mathcal{B}(H)/\mathcal{B}(H)$-measurable. There exist a constant $C<\infty$ such that 
\begin{align*}
&|b(t,y_1,\bar{y}_1)-b(t,y_2,\bar{y}_2)|_{H}+|\sigma(t,y_1,\bar{y}_1)-\sigma(t,y_2,\bar{y}_2)|_{H}+\int_{\textbf{E}}|\theta(t,y_1,\bar{y}_1,e)-\theta(t,y_2,\bar{y}_2,e)|^2\nu(de) \\ 
& \quad  \leq C ( |y_1-y_2|_{H}+|\bar{y}_1-\bar{y}_2|_{H}) \text{ a.s. for all } (\omega,t) \in \Omega \times [0,T].
\end{align*}
We also assume that there exists $C<\infty$ such that:
\begin{align*}
|b(t,y, \bar{y})|^2_{H}+|\sigma(t,y, \bar{y})|^2_{H}+\int_{\textbf{E}}|\theta(t,y,\bar{y},e)|^2_{H} \nu(de) \leq C(1+|y|^2_{H}+|\bar{y}|^2_{H}), \; \forall  (\omega,t) \in  \Omega \times [0,T],  y, \bar{y} \in \textbf{R}.
\end{align*}
Finally, we assume that the operator $\mathbf{F}:\textbf{L}^2(\Omega;H) \mapsto H$ is Fr\'echet differentiable.
% of class $C^1_b$, i.e. $||\nabla \mathbf{F}(\xi)|| \leq C$, for all $\xi \in \textbf{L}_2(\Omega;H)$, where  $||\cdot||$ represents the operatorial norm corresponding to $\mathcal{L}(\textbf{L}_2(\Omega;H);\textbf{L}^*_2(\Omega;H)).$
\end{Assumption}
%Supposing that the above assumption holds, we give the existence and the uniqueness result of a solution in weak sense.
\begin{theorem} Under Assumption \ref{As1}, there exists a unique H-valued progressively measurable process $Y_t, \,\, t \geq 0$ satisfying the mean-field SPDE:
\begin{itemize}
\item[(i)] $Y \in\textbf{L}^2(\Omega \times [0,T], V) \cap \textbf{S}^2(\Omega \times [0,T], H);$
\item[(ii)] $Y_t= h+\int_0^t \left[LY_s+b(s,Y_s,\textbf{F}(Y_s))\right]ds+\int_0^t\sigma(s,Y_s,\textbf{F}(Y_s))dW_s+  \int_0^t\int_{\mathbf{E}} \theta(s_{-},Y_{s_-},\textbf{F}(Y_{s_-}),e) \Tilde{N}(dt,de);$
\item[(iii)] $Y_0=h \in H$.
\end{itemize}
\end{theorem}
\dproof
\textit{\textbf{I. Existence of the solution}}\\
Let $Y_t^0:=h, \,\, t \geq 0.$ For $n \geq 0$, we define $Y^{n+1} \in \textbf{L}^2([0,T]; V) \cap \textbf{S}^2([0,T]; H) $ to be the unique solution to the following equation:
\begin{align}\label{eqqq}
dY_t^{n+1}=LY_t^{n+1}dt+b(t, Y_t^{n+1}, \textbf{F}(Y_t^n))dt+\sigma(t, Y_t^{n+1}, \textbf{F}(Y_t^n))dW_t+\int_{\mathbf{E}} \theta(t_{-}, Y_{t_-}^{n+1}, \textbf{F}(Y_{t_-}^n),e)\Tilde{N}(dt,de).
\end{align}
The solution $Y^{n+1}$ of this equation follows by Proposition 3.1 in \cite{PZ}. Let us now show that the sequence $\{ Y^n, \, n \geq 1 \}$ is a Cauchy sequence in the spaces $\textbf{L}^2(\Omega \times [0,T], V)$ and $\textbf{S}^2(\Omega \times [0,T], H)$. 
By applying It\^{o} formula, we get
\begin{align*}
&|Y_t^{n+1}-Y_t^n|^2_{H} =2 \int_0^t<Y_s^{n+1}-Y_s^n,L(Y_s^{n+1}-Y_s^n)>ds\\  &+2\int_0^t<Y_s^{n+1}-Y_s^n,b(s,Y_s^{n},\mathbf{F}(Y_s^n))-b(s,Y_s^{n-1},\mathbf{F}(Y_s^{n-1}))>_{H}ds \\
&+2\int_0^t<Y_s^{n+1}-Y_s^n,\sigma(s,Y_s^{n+1},\mathbf{F}(Y_s^n))-\sigma(s,Y_s^{n},\mathbf{F}(Y_s^{n-1}))>_{H}dW_s \\
&+\int_0^t|\sigma(s,Y_s^{n},\mathbf{F}(Y_s^n))-\sigma(s,Y_s^{n-1},\mathbf{F}(Y_s^{n-1}))|^2_{H}ds
\\
&+\int_0^t \int_{\textbf{E}}[|\theta(s,Y_{s-}^{n},\mathbf{F}(Y_{s-}^n),e)-\theta(s,Y_{s-}^{n-1},\mathbf{F}(Y_{s-}^{n-1}),e)|^2_{H}]\Tilde{N}(ds,de)
\\
&+2\int_0^t \int_{\textbf{E}}<Y_{s^-}^{n+1}-Y_{s^-}^n,\theta(s,Y_{s^-}^{n},\mathbf{F}(Y_{s^-}^n),e)-\theta(s,Y_{s^-}^{n-1},\mathbf{F}(Y_{s^-}^{n-1}),e)>_{H}\Tilde{N}(ds,de)
\\
&+\int_0^t \int_{\textbf{E}}[|\theta(s,Y_{s-}^{n},\mathbf{F}(Y_{s-}^n),e)-\theta(s,Y_{s-}^{n-1},\mathbf{F}(Y_{s-}^{n-1}),e)|^2_{H} \nu(de)ds.
\end{align*}
Using  Burkholder-Davis-Gundy and Cauchy-Schwarz inequalities and the coercivity assumption \eqref{coerc} on the operator $L$, we obtain that
\begin{align}\label{REL100}
&\mathbf{E}\left[\sup_{0 \leq s \leq t} |Y_s^{n+1}-Y_s^{n}|^2_H\right] \leq -\chi \mathbf{E}[\int_0^t |Y_s^{n+1}-Y_s^n|^2_{V}ds]+C \mathbf{E}[\int_0^t[|Y_s^{n+1}-Y_s^{n}|^2_{H}ds]\nonumber \\
& \quad + \dfrac{1}{2} \mathbf{E}[\sup_{0 \leq s \leq t} |Y_s^{n+1}-Y_s^n|^2_{H}ds]+C \mathbf{E}[\int_0^t[|b(s,Y_s^{n}, \mathbf{F}(Y_s^{n}))-b(s,Y_s^{n-1},\mathbf{F}(Y_s^{n-1}))|^2_{H}ds]\nonumber \\
& \quad +C \mathbf{E}[\int_0^t[|\sigma(s,Y_s^{n}, \mathbf{F}(Y_s^{n}))-\sigma(s,Y_s^{n-1},\mathbf{F}(Y_s^{n-1}))|^2_{H}ds] \nonumber \\
& \quad+ C \mathbf{E}[\int_0^t \int_{\textbf{E}}[|\theta(s,Y_s^{n}, \mathbf{F}(Y_s^{n}),e)-\theta(s,Y_s^{n-1},\mathbf{F}(Y_s^{n-1}),e)|^2_{H}\nu(de)ds].
\end{align}
By the Lipschitz properties of $b, \sigma$ and $\theta$, we deduce
\begin{align}\label{eqqF2}
\mathbf{E}\left[\sup_{0 \leq s \leq t} |Y_s^{n+1}-Y_s^{n}|^2_H\right] \leq & C \mathbf{E}[\int_0^t[|Y_s^{n+1}-Y_s^{n}|^2_{H}ds]+C \mathbf{E}[\int_0^t[|Y_s^{n}-Y_s^{n-1}|^2_{H}ds] \nonumber \\&+C\mathbf{E}[\int_0^t |\mathbf{F}(Y_s^{n})-\mathbf{F}(Y_s^{n-1})|^2_Hds].
\end{align}
We use the mean theorem and obtain the existence for each $n \in \mathbb{N}$, $t \in [0,T]$ of  a random variables $\Tilde{Y}^n(t) \in \textbf{L}^2(\Omega, H)$ such that
\begin{equation}\label{eqF5}
|\mathbf{F}(Y_t^{n})-\mathbf{F}(Y_t^{n-1})|_{H} \leq 
 \|\nabla \mathbf{F}({\Tilde{Y}^n(t)}) \| \|Y_t^n-Y_t^{n-1} \|_{\textbf{L}^2(\Omega;H)}.
\end{equation}
The two above relations \eqref{eqqF2} and \eqref{eqF5}  lead to:
\begin{equation}\label{rel}
\mathbf{E}\left[\sup_{0 \leq s \leq t} |Y_s^{n+1}-Y_s^{n}|^2_H\right] \leq C \mathbf{E}[\int_0^t[|Y_s^{n+1}-Y_s^{n}|^2_{H}ds]+ C \mathbf{E}[\int_0^t |Y_s^{n}-Y_s^{n-1}|_{H}ds].
\end{equation}
Let us now define
$$a^{n}_t=\mathbf{E}\left[ \sup_{0 \leq s \leq t} |Y_s^{n}-Y_s^{n-1}|^2_{H}ds\right]\,\,\,\, A^n_t= \int_0^t a_s^n ds.$$
Using \eqref{rel}, we obtain:
\begin{equation}\label{rel2}
a_t^{n+1} \leq C A_t^{n+1}+CA_t^{n}.
\end{equation}
We multiply the above inequality by $e^{-Ct}$ and derive
$$
\dfrac{d(A_t^{n+1}e^{-Ct})}{dt} \leq Ce^{-Ct}A_t^{n},
$$
which allows us to conclude that
$$A_t^{n+1} \leq Ce^{Ct} \int_0^t e^{-Cs}A_s^n ds \leq Ce^{Ct}tA_t^n.$$
This inequality together with $\eqref{rel2}$ gives
\begin{align*}
a_t^{n+1} \leq C^2 e^{Ct}tA_t^n+CA_t^n \leq C_T \int_0^t A_s^nds,
\end{align*}
where $C_T$ is a given constant. By iteration for all $n$, we finally obtain
$$\mathbf{E}[\sup_{0 \leq s \leq T} |Y_s^{n+1}-Y_s^n|^2_{H}]\leq C \dfrac{(C_T T)^n}{n!}.$$
This implies that we can find $Y \in \textbf{S}_2(\Omega \times [0,T]; H)$ such that
\begin{align*}
\lim_{n \rightarrow \infty} \mathbf{E}[\sup_{0 \leq s \leq t} |Y_s^n-Y_s|_H^2ds]=0.
\end{align*}
By $\eqref{REL100}$, we remark that $Y^n$ also converges to $Y$ in $\textbf{L}^2(\Omega \times [0,T],V).$
Passing to the limit in $\eqref{eqqq},$ we obtain that $Y$ satisfies this equation.

%We now prove the uniqueness of the solution.\\

\textbf{\textit{II. Uniqueness of the solution}}

Let $Y_1$ and $Y_2$ be two solutions in $\textbf{S}^2(\Omega \times [0,T],H) \cap \textbf{L}^2(\Omega \times [0,T],V).$ Applying Ito formula, we have
\begin{align*}
|Y_t^1-Y_t^2|_{H}^2  = &-2 \int_0^t<Y_s^1-Y_s^2,L(Y_s^1-Y_s^2)>ds \\
&+2 \int_0^t<Y_s^1-Y_s^2, b(s,Y_s^1, \textbf{F(}Y_s^1)\textbf{)}-b(s,Y_s^2, \textbf{F(}Y_s^2)\textbf{)}>_Hds  \\
&+2 \int_0^t<Y_s^1-Y_s^2, \sigma(s,Y_s^1, \textbf{F(}Y_s^1)\textbf{)}-\sigma(s,Y_s^2, \textbf{F(}Y_s^2)\textbf{)}>_HdW_s  \\
&+ \int_0^t|\sigma(s,Y_s^1, \textbf{F(}Y_s^1)\textbf{)}-\sigma(s,Y_s^2, \textbf{F(}Y_s^2)\textbf{)}|^2_{H}ds  \\
&+ \int_0^t \int_{\textbf{E}} \left[|\theta(s,Y_{s^-}^1, \textbf{F(}Y_{s^-}^1\textbf{)},e)-\theta(s,Y_{s^-}^2, \textbf{F(}Y_{s^-}^2\textbf{)},e)|^2_{H} \right. \\
&  \left.+2<Y_{s^-}^1-Y_{s^-}^2, \theta(s,Y_{s^-}^1, \textbf{F(}Y_{s^-}^1\textbf{)},e)-\theta(s,Y_{s^-}^2, \textbf{F(}Y_{s^-}^2\textbf{)},e)>\right]\Tilde{N}(ds,de) \\
&+ \int_0^t \int_{\textbf{E}}|\theta(s,Y_{s^-}^1, \textbf{F(}Y_{s^-}^1\textbf{)},e)-\theta(s,Y_{s^-}^2, \textbf{F(}Y_{s^-}^2\textbf{)},e)|^2_{H} ds \nu(de).
\end{align*}
Using now the coercivity assumption on the operator $L$, the Lipschitz property of $b, \sigma, \theta$ and the boundness of the Fr\'{e}chet derivative of the operator $\textbf{F}$, we finally obtain:

\begin{align*}
\mathbf{E}[|Y_t^1-Y_t^2|^2_{H}] \leq -\alpha \mathbf{E}[\int_0^t |Y_s^1-Y_s^2|^2_Vds]+C \mathbf{E}[\int_0^t |Y_s^1-Y_s^2|^2_Hds]\nonumber \\
+\dfrac{1}{2} \mathbf{E}[\sup_{0 \leq s \leq t} |Y_t^1-Y_t^2|^2_H]+C\mathbf{E}[\int_0^t |b(s,Y_s^1,\textbf{F(}Y_s^1)\textbf{)}-b(s,Y_s^2,\textbf{F(}Y_s^2)\textbf{)}|^2_{H}ds]\nonumber \\
+C\mathbf{E}[\int_0^t |\sigma(s,Y_s^1,\textbf{F(}Y_s^1)\textbf{)}-\sigma(s,Y_s^2,\textbf{F(}Y_s^2)\textbf{)}|^2_{H}ds]+ \nonumber \\ C\mathbf{E}[\int_0^t \int_{\textbf{E}} |\theta(s,Y_s^1,\textbf{F(}Y_s^1)\textbf{)}-\theta(s,Y_s^2,\textbf{F(}Y_s^2)\textbf{)}|^2_{H}\nu(de)ds] \\
 \leq C\mathbf{E}[\int_0^t |Y_s^1-Y_s^2|^2_{H}ds].
\end{align*}

We thus deduce that $Y_t^1=Y_t^2.$
\section{Existence and uniqueness of general mean-field backward SPDEs with L\'evy noise}

In this section we give an existence and uniqueness result for mean-field backward SPDEs with jumps.
The analysis will be carried out in a general case, where there exists a \textit{general mean-field operator} acting on each composant of the solution.

 We consider the same framework as in the previous section. Let $A$ be a bounded linear operator from $V$ to $V^*$ satisfying the following coercivity hypothesis: There exist constants $\alpha>0$ and 
$\lambda \geq 0$ such that
$$
2<Au,u>+\lambda|u|^2_{H} \geq \alpha || u||^2_{V}\,\,\, \text{ for all }  u \in V,
$$
where $<Au,u>=Au(u)$ denotes the action of $Au \in V^*$ on $u \in V$.
%Let $K$ be another separable Hilbert space.  Let $\{W_t,\,\,t\geq 0\}$ be a cylindrical Brownian motion  with covariance space $K$ on the probability space $(\Omega, \mathcal{F}, P)$, i.e. for any $k \in K$, $<W_t,k>$ is  a real valued-Brownian motion with $\mathbf{E}[<W_t,k>^2]=t|k|_K^2$. We denote by $L^2(K,H)$ the Hilbert space of Hilbert-Schmidt operators from $K$ into $H$ equipped with the inner product $<S_1,S_2>_{L^2(K,H)}= \sum_{i=1}^{\infty}<S_1k_{i}, S_2k_{i}>_{H}$. Let $(\mathbf{E}, \mathcal{B}(\mathbf{E}))$ be a measurable space, where $\mathbf{E}$ is a topological vector space. Let $\eta(t)$ be a L\'evy process on $\mathbf{E}$. Denote by $\nu(de)$ the L\'evy measure of $\eta$. Denote by $L^2(\nu)$ the $L^2$-space of  square integrable $H$-valued measurable functions associated with $\nu$. Set $p(t)=\eta(t)-\eta(t-)$.\\

\begin{Assumption}\label{As2}
Let  $f:[0,T] \times \Omega \times H \times H \times H \times H \times \textbf{L}^2_{\nu}(H) \times \textbf{L}^2_{\nu}(H) \rightarrow H$ be a $\mathcal{P} \times \mathcal{B}(H) \times \mathcal{B}(H) \times \mathcal{B}(H) \times \mathcal{B}(H) \times \mathcal{B}(\textbf{L}_{\nu}^2(H))\times \mathcal{B}(\textbf{L}_{\nu}^2(H)) / \mathcal{B}(H)$ measurable. 
There exists a constant $C<\infty$ such that
 \begin{align*}
|f(t,\omega,y_1,\Tilde{y}_1, z_1,\Tilde{z}_1, q_1, \Tilde{q}_1)&-f(t, \omega, y_2,\Tilde{y}_2, z_2, \Tilde{z}_2, q_2, \tilde{q}_2)|_{H}  \leq C(|y_1-y_2|_{H}+|\Tilde{y}_1-\Tilde{y}_2|_{H} \\&+|z_1-z_2|_{H} +|\Tilde{z}_1-\Tilde{z}_2|_{H}+|q_1-q_2|_{\textbf{L}^2_{\nu}(H)}+|\Tilde{q}_1-\Tilde{q}_2|_{\textbf{L}^2_{\nu}(H)})
\end{align*}
for all $t,y_1, \Tilde{y}_1,z_1, \Tilde{z}_1, q_1,\Tilde{q}_1,y_2,\Tilde{y}_2,z_2, \Tilde{z}_2, q_2, \Tilde{q}_2.$
We also assume the integrability condition
\begin{align}
\mathbf{E}[\int_0^T |f(t,0,0,0,0,0,0)|^2_{H}dt]<\infty.
\end{align}
\end{Assumption}
We now give our main result of existence and uniqueness.
%Before doing this, we introduce the following three spaces which will be used in the sequel:\\
%
%- $L^2(\Omega,H)$ represents the set of random variables $\mathcal{X}$ such that $\mathbf{E}[|\mathcal{X}|^2_{H}]< \infty;$\\
%
%- $L^2(\Omega,L^2(K,H))$ represents the set of random variables $\mathcal{X}$ such that $\mathbf{E}[|\mathcal{X}|^2_{L^2(K,H)}]< \infty;$\\
%
%- $L^2(\Omega,L^2(\nu))$ represents the set of random variables $\mathcal{X}$ such that $\mathbf{E}[|\mathcal{X}|^2_{L^2(\nu)}]< \infty.$
\begin{theorem}
Assume Assumption \ref{As2} holds.
Let $\xi \in \textbf{L}^2(\Omega;H)$. Let $\mathbf{\mathcal{H}}:\textbf{L}^2(\Omega;H) \mapsto H$, $\mathcal{J}:\textbf{L}^2(\Omega;H) \mapsto H$ and $\mathcal{K}:\textbf{L}^2(\Omega, \textbf{L}^2_\nu(H)) \mapsto \textbf{L}^2_\nu(H)$ be Fr\'echet differentiable operators. There exists a unique $H \times H \times \textbf{L}^2_{\nu}(H)$-valued progressively measurable process $(Y_t,Z_t,U_t)$ such that
$$
\text{(i)}\,\,\,\mathbf{E}[|Y_t|^2_{H}]<\infty, \,\, \mathbf{E}[\int_0^T|Z_t|^2_{H}]<\infty,  \,\, \mathbf{E}[\int_0^T|U_t|^2_{\textbf{L}^2_{\nu}(H)}dt]<\infty.
$$ $$
\text{(ii)}\,\,\, \xi=Y_t+\int_t^T AY_sds+\int_t^Tf(s,Y_s,\mathcal{H}(Y_s),Z_s,\mathcal{J}(Z_s),U_s, \mathcal{K}(U_s))ds+\int_t^T Z_sdW_s+\int_t^T \int_{\textbf{E}} U_s(e) \Tilde{N}(ds,de),
$$
for all $\,\,\, 0 \leq t \leq T$.

The equation $(ii)$ should be understood in the dual space $V^*$.
\end{theorem}
\dproof 
\textbf{\textit{I. Existence of the solution}}\\
Set $Y_t^{0}=0; Z_t^{0}=0; U_t^{0}=0$. We denote by $(Y_t^n,Z_t^n,U_t^n)$ the unique solution of the mean-field backward stochastic equation:
$$
\begin{cases}
dY_t^n &=AY_t^ndt+f(t,Y_t^n, \mathcal{H}(Y_t^{n-1}), Z_t^n, \mathcal{J}(Z_t^{n-1}), U_t^n,\mathcal{K}(U_t^{n-1}))dt+Z_t^ndW_t+\int_{\mathbf{E}} U_t^n(e) \Tilde{N}(dt,de) \\
Y_T^n&=\xi.
\end{cases}$$
The existence and the uniqueness of a solution $(Y_t^n,Z_t^n,U_t^n)$ of such an equation has been proved in \cite{OPZ}. By applying It\^{o}'s formula, we get
\begin{align*}
0&=|Y_T^{n+1}-Y_T^n|^2_{H}\\
&=|Y_t^{n+1}-Y_t^n|^2_{H}+2\int_t^T<A(Y_s^{n+1}-Y_s^{n}),Y_s^{n+1}-Y_s^{n}>ds\\ 
&+2\int_t^T\!\!\!<f(s,Y_s^{n+1}, \mathcal{H}(Y_s^{n}),Z_s^{n+1},\mathcal{J}(Z_s^{n}), U_s^{n+1}, \mathcal{K}(U_s^{n}))\\
&\,\,\,\,\ \ \ \ \ \,\ -f(s,Y_s^{n}, \mathcal{H}(Y_s^{n-1}),Z_s^{n}, \mathcal{J}(Z_s^{n-1}), U_s^{n}, \mathcal{K}(U_s^{n-1})),Y_s^{n+1}-Y_s^n>_{H}ds\\
&+\int_t^T \int_{\mathbf{E}} \left[|Y_{s^-}^{n+1}-Y_{s^-}^n+U_s^{n+1}-U_s^{n}|^2_{H}-|Y_{s^-}^{n+1}-Y_{s^-}^n|^2_{H}\right] \Tilde{N}(ds,de)+\int_t^T \int_{\mathbf{E}} [|U_s^{n+1}(e)-U_s^n(e)|^2_{H}]\nu(de)ds\\ 
&+2 \int_t^T<Y_s^{n+1}-Y_s^n, d(\mathcal{Z}_s^{n+1}-\mathcal{Z}_s^n)>_{H}+\int_t^T |Z_s^{n+1}-Z_s^{n}|^2_{H}ds,
\end{align*}
where $\mathcal{Z}^n_t:=\int_0^t Z^n_s dW_s$.

We thus get, by taking the expectation and using the coercivity assumption on the operator $A$
\begin{align}\label{Myeq}
&\mathbf{E}[|Y_t^{n+1}-Y_t^n|^2_{H}]+\mathbf{E}[\int_t^T |Z_s^{n+1}-Z_s^{n}|^2_{H}ds]+\mathbf{E}[\int_t^T \int_{\mathbf{E}} |U_s^{n+1}-U_s^n|^2_{H}\nu(de)ds]=\\
&-2\mathbf{E}[<A(Y_s^{n+1}-Y_s^{n}),Y_s^{n+1}-Y_s^{n}>ds] \nonumber \\
&-2 \mathbf{E}[\int_t^T<f(s,Y_s^{n+1}, \mathcal{H}(Y_s^{n}),Z_s^{n+1}, \mathcal{J}(Z_s^{n}), U_s^{n+1}, \mathcal{K}(U_s^{n})) \nonumber \\
&\ \ \ \ \ \ \ \ -f(s,Y_s^{n}, \mathcal{H}(Y_s^{n-1}),Z_s^{n},\mathcal{J}(Z_s^{n-1}),U_s^{n}, \mathcal{K}(U_s^{n-1})),Y_s^{n+1}-Y_s^n>ds]\leq \nonumber \\
& \leq \lambda \mathbf{E}[\int_t^T |Y_s^{n+1}-Y_s^{n}|^2_{H}ds]-\alpha \mathbf{E}[\int_t^T |Y_s^{n+1}-Y_s^{n}|^2_{V}ds]- \nonumber \\
& -2 \mathbf{E}[\int_t^T<f(s,Y_s^{n+1}, \mathcal{H}(Y_s^{n}),Z_s^{n+1}, \mathcal{J}(Z_s^{n}), U_s^{n+1}, \mathcal{K}(U_s^{n})) \nonumber \\
&\ \ \ \ \ \ \ \ -f(s,Y_s^{n}, \mathcal{H}(Y_s^{n-1}),Z_s^{n},\mathcal{J}(Z_s^{n-1}),U_s^{n}, \mathcal{K}(U_s^{n-1})),Y_s^{n+1}-Y_s^n>_{H}ds]. \nonumber
\end{align}
%& -2 \mathbf{E}[\int_t^T<f(s,Y_s^{n+1}, \mathcal{H}(Y_s^{n-1}),Z_s^{n+1}, \mathcal{J}(Z_s^{n-1}), U_s^{n+1}, \mathcal{K}(U_s^{n-1}))
%&\ \ \ \ \ \ \ \ -f(s,Y_s^{n}, \mathcal{H}(Y_s^{n-1}),Z_s^{n},\mathcal{J}(Z_s^{n-1}),U_s^{n}, \mathcal{K}(U_s^{n-1})),Y_s^{n+1}-Y_s^n>ds]\\
By using the Cauchy Schwarz inequality and the Lipschitz property of the generator $f$, for each $(t , \omega) \in [0,T] \times \Omega$, we obtain:
\begin{align}\label{eq4}
&<f(s,Y_s^{n+1}, \mathcal{H}(Y_s^{n-1}),Z_s^{n+1}, \mathcal{J}(Z_s^{n-1}), U_s^{n+1}, \mathcal{K}(U_s^{n-1})) \nonumber \\
&-f(s,Y_s^{n}, \mathcal{H}(Y_s^{n-1}),Z_s^{n},\mathcal{J}(Z_s^{n-1}),U_s^{n}, \mathcal{K}(U_s^{n-1})),Y_s^{n+1}-Y_s^n>_{H} \nonumber \\
& \leq |f(s,Y_s^{n+1}, \mathcal{H}(Y_s^{n-1}),Z_s^{n+1}, \mathcal{J}(Z_s^{n-1}), U_s^{n+1}, \mathcal{K}(U_s^{n-1}))-f(s,Y_s^{n}, \mathcal{H}(Y_s^{n-1}),Z_s^{n},\mathcal{J}(Z_s^{n-1}),U_s^{n}, \mathcal{K}(U_s^{n-1}))|_{H} \nonumber \\
&\cdot |Y_s^{n+1}-Y_s^n|_{H} \nonumber \\
& \leq C \left( |\mathcal{H}(Y_s^{n})-\mathcal{H}(Y_s^{n-1})|_{H}+ |\mathcal{J}(Z_s^{n})-\mathcal{J}(Z_s^{n-1})|_{H}+ |\mathcal{K}(U_s^{n})-\mathcal{K}(U_s^{n-1})|_{\textbf{L}^2_{\nu}(H)}\right)|Y_s^{n+1}-Y_s^{n}|_{H} \nonumber \\
& +C \big( |Y_s^{n+1}-Y_s^{n}|_{H}+ |Z_s^{n+1}-Z_s^{n}|_{H}+|U_s^{n+1}-U_s^{n}|_{\textbf{L}^2_{\nu}(H)}\big)|Y_s^{n+1}-Y_s^{n}|_{H}.
\end{align}

We now appeal to  the mean theorem in Hilbert spaces and obtain the existence for each $t \in [0,T]$ of  some random variables $\Tilde{Y}^n(t) \in \textbf{L}^2(\Omega, H)$, $\Tilde{Z}^n(t) \in \textbf{L}^2(\Omega, H),$ $\Tilde{U}^n(t) \in \textbf{L}^2(\Omega, \textbf{L}^2_{\nu}(H))$ such that

\begin{align}\label{eq5}
|\mathcal{H}(Y_t^{n})-\mathcal{H}(Y_t^{n-1})|_{H} \leq 
 \|\nabla \mathcal{H}({\Tilde{Y}^n(t)}) \| \|Y_t^n-Y_t^{n-1} \|_{\textbf{L}^2(\Omega, H)}\nonumber \\
|\mathcal{J}(Z_t^{n})-\mathcal{J}(Z_t^{n-1})|_{H}\leq \|\nabla \mathcal{J}({\Tilde{Z}^n(t)}) \| \|Z_t^n-Z_t^{n-1} \|_{\textbf{L}^2(\Omega, H)} \nonumber \\
|\mathcal{K}(U_t^{n})-\mathcal{K}(U_t^{n-1})|_{H}\leq \|\nabla\mathcal{K}({\Tilde{U}^n(t)}) \| \|U_t^n-U_t^{n-1} \|_{\textbf{L}^2(\Omega,\textbf{L}^2_{\nu}(H))}.
\end{align}

Using $\eqref{Myeq}$, $\eqref{eq4}$, $\eqref{eq5}$ together with the boundness of the Fr\'echet derivatives of the operators $\mathcal{H}, \mathcal{J}, \mathcal{K}$ and the inequality $2ab \leq \varepsilon a^2+ \frac{1}{\varepsilon} b^2$, we obtain:
\begin{align*}
&\mathbf{E}[|Y_t^{n+1}-Y_t^n|^2_{H}]+\mathbf{E}[\int_t^T |Z_s^{n+1}-Z_s^{n}|^2_{H}ds]+\mathbf{E}[\int_t^T \int_{\mathbf{E}} |U_s^{n+1}(e)-U_s^n(e)|^2_{H} \nu(de)ds] \leq \\
& \leq \lambda \mathbf{E}[\int_t^T |Y_s^{n+1}-Y_s^{n}|^2_{H}ds]-\alpha \mathbf{E}[\int_t^T |Y_s^{n+1}-Y_s^{n}|^2_{V}ds]-\\
&+C \varepsilon \mathbf{E}[\int_t^T\left(|Y_s^{n}-Y_s^{n-1}|^2_H+ |Z_s^{n}-Z_s^{n-1}|^2_{H}+|U_s^{n}-U_s^{n-1}|^2_{\textbf{L}^2_{\nu}(H)}\right)ds]+\frac{1}{\varepsilon}\mathbf{E}[\int_t^T|Y_s^{n+1}-Y_s^{n}|^2_{H}ds]\\
& +C \beta \mathbf{E}[\int_t^T \left( |Y_s^{n+1}-Y_s^{n}|^2_H+|Z_s^{n+1}-Z_s^{n}|^2_{H}+|U_s^{n+1}-U_s^{n}|^2_{\textbf{L}^2_{\nu}(H)}\right) ds]+\frac{1}{\beta}\mathbf{E}[\int_t^T|Y_s^{n+1}-Y_s^{n}|^2_{H}ds],
\end{align*}
where $C$ is a constant dependind on the Lipschitz constant of $f$ and the bounding constants of the Fr\'echet derivative operators of $\mathcal{H},\mathcal{J},\mathcal{K}$.\\
Let us choose $\varepsilon \leq \frac{1}{4C}$ and $\beta \leq \frac{1}{2C}$. We set $\gamma:=\lambda+C \beta+\frac{1}{\varepsilon}+\frac{1}{\beta}+\frac{1}{2}$ and then multiply the previous inequality by $e^{\gamma t}$. We thus get
\begin{align}\label{eq2}
&-\frac{d}{dt} \left(e^{\gamma t} \mathbf{E}[\int_t^T|Y_s^{n+1}-Y_s^{n}|^2_{H}ds]  \right)+\frac{1}{2} e^{\gamma t} \mathbf{E}[\int_t^T |Z_s^{n+1}-Z_s^{n}|^2_{H}ds] \nonumber \\
&+\frac{1}{2} \mathbf{E}[\int_t^T |Y_s^{n+1}-Y_s^{n}|^2_{H}ds] e^{\gamma t}+\frac{1}{2} e^{\gamma t} \mathbf{E}[\int_t^T  |U_s^{n+1}-U_s^{n}|^2_{\textbf{L}^2_{\nu}(H)}ds]+\alpha e^{\gamma t} \mathbf{E}[\int_t^T |Y_s^{n+1}-Y_s^{n}|^2_{V}ds] \nonumber \\
& \leq \frac{1}{4} \mathbf{E}[\int_t^T |Y_s^{n}-Y_s^{n-1}|^2_{H}] e^{\gamma t}+\frac{1}{4} \mathbf{E}[\int_t^T |Z_s^{n}-Z_s^{n-1}|^2_{H}] e^{\gamma t}+\frac{1}{4} \mathbf{E}[\int_t^T  |U_s^{n}-U_s^{n-1}|^2_{\textbf{L}^2_{\nu}(H)}ds] e^{\gamma t}.
\end{align}
We now integrate between $0$ and $T$ and obtain:
\begin{align}\label{eq1}
 \mathbf{E}[\int_0^T |Y_s^{n+1}-Y_s^{n}|^2_{H}ds]+\frac{1}{2}\int_0^T \mathbf{E}[\int_t^T |Y_s^{n+1}-Y_s^{n}|^2_{H}ds] e^{\gamma t} dt+\frac{1}{2}\int_0^T \mathbf{E}[\int_t^T |Z_s^{n+1}-Z_s^{n}|^2_{H}ds] e^{\gamma t} dt\nonumber \\
+\frac{1}{2} \int_0^T e^{\gamma t} \mathbf{E}[\int_t^T|U_s^{n+1}-U_s^{n}|^2_{\textbf{L}^2_{\nu}(H)}ds]+\int_0^T \alpha \mathbf{E}[\int_t^T |Y_s^{n+1}-Y_s^{n}|^2_{V}ds]e^{\gamma t} dt \nonumber \\
 \leq \frac{1}{4} \int_0^T \mathbf{E}[\int_t^T |Y_s^{n}-Y_s^{n-1}|^2_{H}ds] e^{\gamma t}dt+\frac{1}{4} \int_0^T \mathbf{E}[\int_t^T |Z_s^{n}-Z_s^{n-1}|^2_{H}ds]e^{\gamma t}dt \nonumber \\ + \frac{1}{4} \int_0^T \mathbf{E}[\int_t^T  |U_s^{n}-U_s^{n-1}|^2_{\textbf{L}^2_{\nu}(H)}ds] e^{\gamma t}dt.
\end{align}
From the above inequality it follows that
\begin{align*}
\int_0^T \mathbf{E}[\int_t^T |Y_s^{n}-Y_s^{n-1}|^2_{H}] e^{\gamma t}dt&+ \int_0^T \mathbf{E}[\int_t^T |Z_s^{n}-Z_s^{n-1}|^2_{H}]e^{\gamma t}dt+\int_0^T \mathbf{E}[\int_t^T  |U_s^{n}-U_s^{n-1}|^2_{\textbf{L}^2_{\nu}(H)}ds] e^{\gamma t}dt  \\
& \leq \frac{1}{2^{n}}C.
\end{align*}
From \eqref{eq1} one can deduce
$$
\mathbf{E}[\int_0^T |Y_s^{n+1}-Y_s^{n}|^2_{H}ds] \leq \frac{1}{2^{n}}C.
$$
We now appeal to $\eqref{eq2}$ and derive
\begin{align*}
&\frac{1}{2}\mathbf{E}[\int_0^T |Y_s^{n+1}-Y_s^{n}|^2_{H}ds]+ \frac{1}{2}\int_0^T \mathbf{E}[\int_0^T |Z_s^{n+1}-Z_s^{n}|^2_{H}ds]+\frac{1}{2}\int_0^T \mathbf{E}[\int_0^T  |U_s^{n}-U_s^{n-1}|^2_{\textbf{L}^2_{\nu}(H)}ds]\\
& \leq \gamma  \frac{1}{2^{n}}C+\frac{1}{4}\mathbf{E}[\int_0^T |Y_s^{n}-Y_s^{n-1}|^2_{H}ds]+ \frac{1}{4}\mathbf{E}[\int_0^T |Z_s^{n}-Z_s^{n-1}|^2_{H}ds]+\frac{1}{4} \mathbf{E}[\int_0^T  |U_s^{n+1}-U_s^{n}|^2_{\textbf{L}^2_{\nu}(H)}ds],
\end{align*}
which implies that
$$ \mathbf{E}[\int_0^T |Y_s^{n+1}-Y_s^{n}|^2_{H}ds]+\int_0^T \mathbf{E}[\int_0^T |Z_s^{n+1}-Z_s^{n}|^2_{H}ds]+\int_0^T \mathbf{E}[\int_0^T  |U_s^{n}-U_s^{n-1}|^2_{\textbf{L}^2_{\nu}(H)}ds] \leq \frac{1}{2^{n-1}}C \gamma n. $$
This leads to
$$
\mathbf{E}[\int_0^T |Y_s^{n+1}-Y_s^n|^2_Vds] \leq (\frac{1}{2})^{n-1}(n+1)C\gamma.
$$ Hence, we can conclude that the sequence $(Y^n,Z^n,U^n)$, ${n \geq 1}$  is a Cauchy sequence in the Banach space $L^2(\Omega \times [0,T], V) \times L^2(\Omega \times [0,T], H) \times L^2(\Omega \times [0,T], L^2(\nu))$, and thus converges in the corresponding spaces to $(Y,Z,U)$. The limit $(Y,Z,U)$ satisfies:
$$
Y_t+\int_t^TAY_sds+\int_t^Tf(s,Y_s,\mathcal{H}(Y_s), Z_s,\mathcal{J}(Z_s),U_s,\mathcal{K}(U_s))ds+\int_t^TZ_sdW_s+\int_t^T \int_\mathbf{E} U_s\Tilde{N}(ds,de)=\xi\,\, \text{ a.s. }
$$

\textbf{\textit{II. Uniqueness of the solution}}\\
The proof of the uniqueness of the solution is classical, but we give it for convenience of the reader. Suppose  $(Y_t,Z_t,U_t)$ and $(\Tilde{Y}_t,\Tilde{Z}_t,\Tilde{U}_t)$ are two solutions.
By applying It\^{o} formula, we obtain
\begin{align*}
&\mathbf{E}[|Y_t-\Tilde{Y}_t|^2_{H}]+\mathbf{E}[\int_t^T |Z_s-\Tilde{Z}_s|^2_{H}ds]+\mathbf{E}[\int_t^T |U_s-\Tilde{U}_s|^2_{\textbf{L}^2_{\nu}(H)}]ds]=\nonumber \\
&-\mathbf{E}[<A(Y_s-\Tilde{Y}_s),Y_s-\Tilde{Y}_s>ds] \nonumber \\
&-2 \mathbf{E}[\int_t^T<f(s,Y_s, \mathcal{H}(Y_s),Z_s, \mathcal{J}(Z_s), U_s, \mathcal{K}(U_s)) \nonumber \\
&\ \ \ \ \ \ \ \ -f(s,\Tilde{Y}_s, \mathcal{H}(\Tilde{Y}_s), \Tilde{Z}_s,\mathcal{J}(\Tilde{Z}_s),\Tilde{U}_s, \mathcal{K}(\Tilde{U}_s)),\Tilde{Y}_s-\Tilde{Y}_s>_{H}ds] \nonumber \\
& \leq \lambda \mathbf{E}[\int_t^T |Y_s-\Tilde{Y}_s|^2_{H}ds]-\alpha \mathbf{E}[\int_t^T |Y_s-\Tilde{Y}_s|^2_{V}ds]+K\mathbf{E}[\int_t^T|Y_s-\Tilde{Y}_s|^2_{H}ds] \nonumber \\
& +\frac{1}{2}\mathbf{E}[\int_t^T|Z_s-\Tilde{Z}_s|^2_{H}ds]+\frac{1}{2}\mathbf{E}[\int_t^T|U_s-\Tilde{U}_s|^2_{\textbf{L}^2_{\nu}(H)}ds].
\end{align*}
We thus derive that
$$
\mathbf{E}[|Y_t-\Tilde{Y}_t|_{H}^2] \leq (\lambda+K) \mathbf{E}[\int_t^T|Y_s-\Tilde{Y}_s|_{H}^2]
$$
Hence, by Gronwall lemma, we get
$$Y_t=\Tilde{Y}_t.$$
This also implies that $Z_t=\Tilde{Z}_t$ and $U_t=\Tilde{U}_t.$
\appendix
\section{Some results on Banach theory}

We recall here some basic concepts and results from Banach space theory. Let $V$ be an open subset of a Banach space $\mathcal{X}$ with norm $ \|\cdot \|$ and let $F:V\mapsto \mathbb{R}$.

\textbf{(i)} We say that $F$ has a directional derivative (or G\^{a}teaux derivative ) at $x \in \mathcal{X}$ in the direction $y \in \mathcal{X}$ if
$$
D_{y}F(x):= \lim_{\varepsilon \rightarrow 0}\dfrac{1}{\varepsilon}(F(x+\varepsilon y)-F(x))
$$
exists.\\
\textbf{(ii)} We say that $F$ is Fr\'echet differentiable at $x \in V$ if there exists a linear map $L:\mathcal{X}\mapsto \mathbb{R}$ such that
$$
\lim_{h \rightarrow 0; h \in \mathcal{X}} \dfrac{1}{ \| h \|}|F(x+h)-F(x)-L(h)|=0.
$$
In this case we call $L$ the \textit{gradient} (or Fr{\'echet derivative}) of $F$ at $x$ and we write 
$$
L=\nabla F.
$$
\textbf{(iii)} If $F$ is Fr{\'e}chet differentiable, then $F$ has a directional derivative in all directions $y \in \mathcal{X}$ and
$$
D_yF(x)= \nabla_xF(y)=:<\nabla_xF,y>.
$$ In particular, if $\mathcal{X}=\textbf{L}^2(P)$ the Fr\'echet derivative of $F$ at $X \in \textbf{L}^2(P)$, denoted by $\nabla F(X)$, is a bounded linear functional on $\textbf{L}^2(P)$, which we can identify by Riesz theorem with a random variable in $\textbf{L}^2(P)$. For example, if $F(X)=\mathbf{E}[\phi(X)]; \,\, X \in \textbf{L}^2(P),$ where $\phi$ is a real $C^1$-function such that $\phi(X)\in \textbf{L}^2(P)$ and $\dfrac{\partial \phi}{\partial x}(X) \in L^2(P)$, then $\nabla F(X)= \dfrac{\partial \phi}{\partial x}(X)$ and \\
$\nabla F(X)(Y)=\left<\dfrac{\partial \phi}{\partial x}(X),Y\right>_{\textbf{L}^2(P)}=\mathbf{E}\left[\dfrac{\partial \phi}{\partial x}(X)Y\right],$ for $Y \in \textbf{L}^2(P).$

%It remains to:
%- to put the condition on $b(0,0,0)$\\
%- to put the controls square integrable\\
%- to put $U$ convex\\
%- the integrability condition to remove\\
%- example for application of the results  of the existence and uniqueness to the first part (Sobolev space with the boundary conditions)

%\fproof

\end{document}